\documentclass{amsart}

\usepackage{amsfonts}
\usepackage{amsmath}
\usepackage{amsthm}
\usepackage{amssymb}
\usepackage[english]{babel}
\usepackage{graphics}
\usepackage{latexsym}
\usepackage{longtable} \usepackage{mathrsfs}
\usepackage{graphicx}
\usepackage{wrapfig} 
\usepackage[pdftex,colorlinks=true,linkcolor=blue,citecolor=cyan]{hyperref}

\newtheorem{theorem}{Theorem}
\newtheorem{corollary}[theorem]{Corollary}
\newtheorem{lemma}[theorem]{Lemma}

\newtheorem{claim}[theorem]{Claim}
\newtheorem{example}[theorem]{Example}

\theoremstyle{definition}
\newtheorem{definition}[theorem]{Definition}
\newtheorem{remark}[theorem]{Remark}

\newcommand{\msL}{\mathscr{L}}

\newcommand{\mF}{\mathcal{F}}

\newcommand{\mM}{\mathcal{M}}

\newcommand{\mY}{\mathscr{Y}}
\newcommand{\mD}{\mathcal{D}}

\newcommand{\R}{\mathbb{R}}
\newcommand{\N}{\mathbb{N}}

\renewcommand{\P}{\mathrm{P}}

\newcommand{\D}{\mathrm{D}}

\newcommand{\noi}{\noindent}
\newcommand{\ms}{\medskip}

\newcommand{\al}{\alpha}
\newcommand{\be}{\beta}
\newcommand{\ga}{\gamma}
\newcommand{\de}{\delta}
\newcommand{\De}{\Delta}
\newcommand{\e}{\varepsilon}
\newcommand{\si}{\sigma}
\newcommand{\la}{\lambda}

\newcommand{\Om}{\Omega}
\newcommand{\om}{\omega}

\newcommand{\weak }{\, -\!\!\!\!-\!\!\!\rightharpoonup}
\newcommand{\weakstar }{ \overset{\, *_{\phantom{|}}}{{\smash{\weak }}\, } }

\newcommand{\larrow}{\longrightarrow}
\newcommand{\ot}{\otimes}

\newcommand{\lmapsto}{\longmapsto}
\newcommand{\ri}{\rightarrow}

\newcommand{\p}{\partial}
\newcommand{\sub}{\subseteq}
\newcommand{\set}{\setminus}
\newcommand{\by}{\times}

\newcommand{\sgn}{\mathrm{sgn}}
\newcommand{\ess}{\mathrm{ess}}

\newcommand{\inter}{\mathrm{int}}

\newcommand{\supp}{\mathrm{supp}}

\newcommand{\bt}{\begin{theorem}}\newcommand{\et}{\end{theorem}}
\newcommand{\bd}{\begin{definition}}\newcommand{\ed}{\end{definition}}
\newcommand{\bl}{\begin{lemma}}\newcommand{\el}{\end{lemma}}
\newcommand{\beq}{\begin{equation}}\newcommand{\eeq}{\end{equation}}
\newcommand{\bc}{\begin{claim}}\newcommand{\ec}{\end{claim}}
\newcommand{\bex}{\begin{example}}\newcommand{\eex}{\end{example}}
\newcommand{\bcor}{\begin{corollary}}\newcommand{\ecor}{\end{corollary}}
\newcommand{\bp}{\begin{proof}}\newcommand{\ep}{\end{proof}}

\newcommand{\BPL}{\medskip \noindent \textbf{Proof of Lemma} }

\numberwithin{equation}{section}
\numberwithin{theorem}{section}

\begin{document}

\title[Generalised solutions and Vectorial Calculus of Variations in $L^\infty$]{Absolutely Minimising Generalised solutions to the Equations of Vectorial Calculus of Variations in $L^\infty$}

\author{Nikos Katzourakis}
\address{Department of Mathematics and Statistics, University of Reading, Whiteknights, PO Box 220, Reading RG6 6AX, Reading, UK}
\email{n.katzourakis@reading.ac.uk}

    \thanks{\!\!\!\!\!\!\texttt{The author has been partially financially supported by the EPSRC grant EP/N017412/1}}


\date{}



\begin{abstract} Consider the supremal functional
\[ \tag{1} \label{1}
E_\infty(u,A) \,:=\, \|\msL(\cdot,u,\D u)\|_{L^\infty(A)},\quad A\subseteq \Omega,
\]
applied to $W^{1,\infty}$ maps $u:\Om\sub \R\longrightarrow \R^N$, $N\geq 1$. Under certain assumptions on $\msL$, we prove for any given boundary data the existence of a map which is: 

\noindent i) a vectorial Absolute Minimiser of \eqref{1} in the sense of Aronsson, 

\noindent ii) a generalised solution to the ODE system associated to  \eqref{1} as the analogue of the Euler-Lagrange equations, 

\noindent iii) a limit of minimisers of the respective $L^p$ functionals as $p\rightarrow\infty$ for any $q\geq 1$ in the strong $W^{1,q}$ topology \&

\noindent iv) partially $C^2$ on $\Om$ off an exceptional compact nowhere dense set. 

\noi {Our method is based on $L^p$ approximations and stable a priori partial regularity estimates. For item ii) we utilise the recently proposed by the author notion of $\mathcal{D}$-solutions in order to characterise the limit as a generalised solution. Our results are motivated from and apply to Data Assimilation in Meteorology.}
\end{abstract}

\maketitle

\tableofcontents

\section{Introduction} \label{section1}

Calculus of Variations in $L^\infty$ has a long history and was pioneered by Aronsson in the 1960s \cite{A1}-\cite{A5}. In the vector case and in one spatial dimension, the basic object of study is the functional
\beq \label{1.1}
E_\infty (u,A)\, :=\, \big\| \msL(\cdot,u,\D u)\big\|_{L^\infty(A)}, \quad u\, :\, \Om\sub \R \larrow \R^N,\ A\sub \Om.
\eeq
Here $u\in W^{1,\infty}_{\text{loc}}(\Om,\R^N)$, $N\geq1$, $\Om$ is an open interval, $A$ is measurable and $\msL \in C^2\big(\Om \by \R^N \by  \R^N\big)$ is a function which we call Lagrangian and whose arguments will be denoted by $(x,\eta,P)$. Aronsson who studied the case $N=1$ was the first to note the locality problems associated to this functional. By introducing the appropriate minimality notion in $L^\infty$, among other things proved the equivalence between the so-called Absolute Minimisers and solutions of the analogue of the Euler-Lagrange equation which is associated to the functional under $C^2$ smoothness hypotheses. The minimality notion of Aronsson adapted to the vector case of \eqref{1.1} is
\beq \label{1.1a}
\ \ \ E_\infty (u,\Om')\,\leq\, E_\infty (u+\phi,\Om'),\quad \forall \,\Om'\Subset \Om,\,\phi\in W^{1,\infty}_0(\Om',\R^N).
\eeq
The higher dimensional scalar analogue when $u :\Om\sub \R^n \ri\R$ is a real function has also attracted considerable attention by the community and by now there is a vast literature, for instance see Crandall \cite{C}, Barron-Evans-Jensen \cite{BEJ}, and for a pedagogical introduction see \cite{K0} and references therein. In particular, the Crandall-Ishii-Lions theory of Viscosity solutions proved to be an indispensable tool in order to study the equations in $L^\infty$ which are non-divergence, highly nonlinear and degenerate. Even in the simplest case where the Lagrangian is the Euclidean norm, i.e.\ $\msL(P)=|P|^2$, in general the solutions are non-smooth and the corresponding PDE which is called $\infty$-Laplacian for smooth functions reads
\beq \label{1.2}
\De_\infty u\, :=\, \D u \ot \D u :\D^2u\, =\sum_{i,j=1}^n\D_i u\, \D_ju\,\D^2_{ij}u\,=\, 0.
\eeq

{However, until the early 2010s, the theory was essentially restricted to the scalar case $N=1$. A most notable exception is the early vectorial contributions of Barron-Jensen-Wang \cite{BJW1, BJW2}.} Therein the authors among other far-reaching results studied the weak* lower semicontinuity of general supremal functionals and proved under certain assumptions the existence of Absolute Minimisers in the ``rank-1" cases, i.e.\ when either $n=1$ or $N=1$.

{In a series of recent papers \cite{K1}-\cite{K9}, the author has initiated the systematic study of the vector-valued case, which except for its intrinsic mathematical interest, appears to be important for many real-world applications (see also the joint contributions with Abugirda, Croce, Pisante and Pryer \cite{AK, CKP, KP, KP2}). In particular, the complete PDE system arising in $L^\infty$ was derived and studied in \cite{K1}}. The results in the aforementioned papers include in particular the study of the analytic properties of classical solutions to the fundamental equations and their connection to the supremal functional. In the case of 
\beq  \label{1.3}
\ \ \ \ \ \ E_\infty(u,A)\, =\, \big\| |\D u|^2 \big\|_{L^\infty(A)},\quad \ u \,:\ \Om \sub \R^n \larrow \R^N,\ \, A\sub \Om,
\eeq
(where $|\D u|$ denotes the Euclidean norm of the gradient on $\R^{N\by n}$), the respective $\infty$-Laplace system written for smooth maps is
\beq \label{1.4}
\De_\infty u \, :=\, \Big(\D u \ot \D u\, +\, |\D u|^2[\D u]^\bot \! \ot I \Big):\D^2u\, =\, 0.
\eeq
In \eqref{1.4}, $[\D u(x)]^\bot$ denotes the orthogonal projection on the orthogonal complement of the range of the linear map $\D u(x) : \R^n \larrow \R^N$ and in index form reads
\[
\begin{split}
\sum_{\be=1}^N\sum_{i,j=1}^n \Big(\D_i u_\al \, \D_ju_\be\ +&\ |\D u|^2[\D u]_{\al \be}^\bot \, \de_{ij}\Big)\, \D_{ij}^2u_\be\, =\, 0, \ \ \ \al=1,...,N,
\\
&[\D u]^\bot\,:=\, \textrm{Proj}_{R(\D u)^\bot}.
\end{split}
\]
An extra difficulty of \eqref{1.4} which is not present in the scalar case of \eqref{1.2} is that the coefficients may be discontinuous along interfaces even for $C^\infty$ solutions because the term involving $[\D u]^\bot$ measures the dimension of the tangent space of $u(\Om)\sub \R^N$ (see \cite{K1, K2} and the numerical experiments in \cite{KP}). This is a general vectorial phenomenon studied in some detail in \cite{K3}. The appropriate minimality notion allowing to connect \eqref{1.4} to the functional \eqref{1.3} has been established in \cite{K4}. It is a remarkable fact that \emph{when the rank of the gradient is greater than one, i.e.\ for maps $u:\Om\sub \R^n \larrow \R^N$ such that $\min\{n,N\}\geq 2$, {Absolute Minimimality for \eqref{1.3} is neither necessary nor sufficient for solvability of \eqref{1.4}} and the correct notion of $\infty$-Minimal maps is intrinsically different} (see \cite{K4}). In the case of the supremal functional \eqref{1.1}, the associated equations written for smooth maps $u:\Om\sub \R \larrow \R^N$ read
\beq \label{1.5a}
\mathcal{F}_\infty \big(\cdot,u,\D u,\D^2u\big)\,=\,0, \quad \text{on }\Om,
\eeq
where
\beq  \label{1.5b}
\begin{split}
\mathcal{F}_\infty (x,\eta,P,X) \,  :=&  \, \Big[  \msL_P(x,\eta,P) \ot \msL_P(x,\eta,P)\\
 &  +  \, \msL(x,\eta,P) [\msL_P(x,\eta,P)]^\bot  \msL_{PP}(x,\eta,P)\Big]X  
\\ 
 &+\,  \Big( \msL_\eta(x,\eta,P) \cdot P \, +\, \msL_x(x,\eta,P) \Big) \msL_P  (x,\eta,P)\\
  & + \,  \msL  (x,\eta,P)\big[\msL_P(x,\eta,P)\big]^\bot \Big(\msL_{P\eta}(x,\eta,P)P
 \\
 & +  \,  \msL_{Px}  (x,\eta,P) \,-\, \msL_\eta(x,\eta,P) \Big).
\end{split}
\eeq
In \eqref{1.5b}, the notation of subscripts denotes derivatives with respect to the respective variables and $\big[\msL_P(x,\eta,P)\big]^\bot$ is the orthogonal projection
\beq  \label{1.5c}
\ \ \big[\msL_P(x,\eta,P)\big]^\bot\, :=\, I - \, \sgn\big(\msL_P(x,\eta,P)\big) \ot \sgn\big(\msL_P(x,\eta,P)\big).
\eeq
The system \eqref{1.5a}-\eqref{1.5c} is a 2nd order ODE system which is quasilinear, non-divergence, non-monotone and with discontinuous coefficients. Even in the scalar case of $N=1$ in which $\mathcal{F}_\infty$ simplifies to
\[
\begin{split}
\mathcal{F}_\infty (x,\eta,P,X) \,  =\,  \big(\msL_P(x,\eta,P)\big)^2 X  
\, +\, 
\Big( \msL_\eta(x,\eta,P)P + \msL_x(x,\eta,P) \Big) \msL_P  (x,\eta,P)
\end{split}
\] 
it is known since the work of Aronsson that in general does not have solutions any more regular than at best $C^1(\Om,\R^N)$ and their ``weak" interpretation is an issue. {Let us also note that, also inspired by Aronsson's work, Sheffield-Smart \cite{SS} made a vectorial breakthrough relevant to \eqref{1.4} and \eqref{1.3} which was simultaneous to \cite{K1}. They studied smooth vector-valued optimal Lipschitz extensions of functions, deriving a different more singular version of $\infty$-Laplacian than \eqref{1.4}, corresponding to \eqref{1.3} but when the matrix norm of $\D u$ is the nonsmooth operator norm on $\R^{N\by n}$.}

In this paper we study the functional \eqref{1.1}, the associated nonlinear system \eqref{1.5a}-\eqref{1.5c} and their connection. Our main result establishes for any given endpoint data on $\Om$ the existence of a vectorial Absolute Minimiser $u^\infty \in W^{1,\infty}_b(\Om,\R^N)$ of \eqref{1.1} which also is a generalised solution to \eqref{1.5a}-\eqref{1.5c} in a certain new sense to be made precise below. We moreover glean extra information about $u^\infty$; it a partially $C^2$ map off a singular compact set and also is a limit of minimisers of the respective $L^m$ functionals
\beq \label{1.9}
\ \ \ \ \ E_m (u,\Om)\, :=\, \int_\Om \msL(\cdot,u,\D u)^m, \quad u\, :\, \Om\sub \R \larrow \R^N,
\eeq
in the strong $W^{1,q}$ topology as $m\ri \infty$ for any $q\geq 1$. Our results have been motivated from and apply to variational Data Assimilation in the form used in the Earth Sciences and in particular Weather Forecasting. Below we discuss the essential idea of our new notion of generalised solution, our assumptions, state our main result and also draw connections to Data Assimilation.

Motivated in part by the nonlinear systems arising in $L^\infty$, in the very recent paper \cite{K8} the author proposed a new theory of generalised solutions which applies to fully nonlinear PDE systems. In addition, this theory allows to interpret merely measurable general mappings $u:\Om\sub \R^n\larrow \R^N$ as solutions of PDE systems which may even be defined by discontinuous nonlinearities and can be of any order. Our approach is duality-free and bypasses the insufficiency of the standard duality ideas to apply to even linear non-divergence equations with rough coefficients. The standing idea of the use of integration-by-parts in order to pass derivatives to test functions is replaced by a probabilistic description of the limiting behaviour of the difference quotients. This builds on the use of \emph{Young measures valued into compact tori}, which is the compactification of the space wherein the derivatives are valued. Background material on Young measures can be found e.g.\ in \cite{FG, CFV, V, FL, P, E, M}, but for the convenience of the reader we recall herein the rudimentary properties we actually utilise.

The essential idea of our new notion of solution for the case needed in this paper can be briefly motivated as follows. Assume that $u : \Om \sub \R \larrow \R^N$ is a strong a.e.\ solution of the system
\beq  \label{1.8}
\mathcal{F}\big(\cdot,u,\D u,\D^2 u\big)\,=\, 0, \quad \text{on }\Om,
\eeq
in $W^{2,\infty}_{\text{loc}}(\Om,\R^N)$. We need a notion of solution which makes sense even if $u$ is merely $W^{1,\infty}_{\text{loc}}(\Om,\R^N)$. To this end we rewrite the above statement that $u$ is a strong solution in the following unconventional fashion
\beq \label{1.11}
\sup_{X\in \supp(\de_{\D^2 u(x)})}\Big|\mathcal{F}\Big(x,u(x),\D u(x),X\Big)\Big|\,=\, 0, \quad \text{a.e. }x\in \Om.
\eeq
That is, we change from the classical viewpoint that the $2$nd derivative is a map $\D^2u : \Om\sub\R \larrow \R^N$ valued in $\R^N$ to that it is a probability valued map given by the Dirac mass at $\D^2u$:
\[
\de_{\D^2u}\ : \ \Om\sub\R \larrow \mathscr{P}(\R^N),\quad    x\lmapsto \de_{\D^2u(x)}.
\]
Obviously, ``supp" denotes the support of the probability measure. Similarly, if $\D^{1,h}$ stands for the difference quotient operator, it can be shown that we may rewrite the definition of $\D^2u$ as
\beq \label{1.12}
\de_{\D^{1,h}\D u} \weakstar\, \de_{\D^2u}, \quad \text{as }h\ri0.
\eeq
The weak* convergence above is meant in the so-called Young measures into $\R^N$, that is the probability-valued maps $\vartheta : \Om\sub\R \larrow \mathscr{P}(\R^N)$ which are weakly* measurable (see the next section for details and Lemma \ref{lemma2a}). The rationale behind the reformulation \eqref{1.11}-\eqref{1.12} is that {\emph{we may allow more general probability-valued maps arising as ``diffuse" $2$nd derivatives for maps classically differentiable only once.} The latter of course may no longer be the concentration measure-valued maps $\de_{\D^2u}$.} This is indeed possible if we augment $\R^N$ and replace it by its $1$-point compactification $\R^N\cup\{\infty\}$. By considering \emph{Young measures valued into spheres}, we obtain the necessary compactness and the maps $(\de_{\D^{1,h}\D u})_{h\neq 0}$ \textbf{always} have subsequential weak* limits
\[
\de_{\D^{1,h_i}\D u} \weakstar\, \mD^2u, \quad \text{as }h_i\ri0
\]
in the space of sphere-valued Young measures $\Om\sub\R \larrow \mathscr{P}\big(\R^N\cup\{\infty\}\big)$ (even when $u$ is merely once differentiable). The above ideas motivate the notion of \textbf{$\mD$-solutions} and \textbf{diffuse derivatives} in the special case of $W^{1,\infty}$ solutions to \eqref{1.8} (see Definitions \ref{definition7} and \ref{definition11}) and will be taken as principal in this work.

As a first application of this new approach, in the paper \cite{K8} among other things we proved existence of $\mD$-solutions to the Dirichlet problem for \eqref{1.4} when $n=N$. Further results in the context of $\mD$-solutions have been established in \cite{K9,K10,K11}, including certain uniqueness assertions. Herein we focus on \eqref{1.1} and \eqref{1.5a}-\eqref{1.5c}. This is a non-trivial task even in the $1D$ case. In fact, it is not possible to work in the generality of \eqref{1.1}, \eqref{1.5a}-\eqref{1.5c} without structural conditions on $\msL$. The most important restriction is that the Lagrangian has to be radial in $P$. This means that $\msL$ can be written as
\beq \label{1.13}
\msL(x,\eta,P)\, =\, \mathscr{H}\left(x,\eta, \frac{1}{2}|P-\mathscr{V}(x,\eta)|^2 \right).
\eeq
for some mappings $\mathscr{H} : \Om \by \R^N \by [0,\infty) \larrow \R$ and $\mathscr{V}:  \Om \by \R^N   \larrow \R^N$. This condition is justified by the results of \cite{K3} since as we proved therein it is both necessary and sufficient for the ODE system to be degenerate elliptic. The extra bonus is that then the \emph{coefficients of \eqref{1.5b} match and become continuous}. This is a special occurrence due to the $1D$ nature of the problem and can not happen when $n\geq 2$. Anyhow, under the assumption \eqref{1.13}, \eqref{1.5b} after a certain rescaling becomes
\begin{align}  \label{1.14}
\mF_\infty(x,\eta,P,X)  & \, = \, \big|P-\mathscr{V} (x,\eta)\big|^2 \mathscr{H}_p \left(x,\eta, \frac{1}{2}\big|P-\mathscr{V}(x,\eta)\big|^2 \right)^{\!\!2} \Big[ X -\mathscr{V}_\eta(x,\eta) P
 \nonumber\\
&  -\mathscr{V}_x(x,\eta) \Big] \, + \, \big(P-\mathscr{V} (x,\eta)\big) \mathscr{H}_p\left(x,\eta, \frac{1}{2}\big|P-\mathscr{V}(x,\eta)\big|^2 \right)  
\nonumber \\
& \centerdot  \Bigg[\mathscr{H}_x\left(x,\eta, \frac{1}{2}\big|P-\mathscr{V}(x,\eta)\big|^2 \right)
\, +\, P\cdot \mathscr{H}_\eta\left(x,\eta, \frac{1}{2}\big|P-\mathscr{V}(x,\eta)\big|^2 \right) \Bigg]
 \\
  & -\,  \mathscr{H}_p\left(x,\eta, \frac{1}{2}\big|P-\mathscr{V}(x,\eta)\big|^2 \right) \Bigg( \big|P-\mathscr{V} (x,\eta)\big|^2I
  \nonumber \\ 
&- \big(P-\mathscr{V} (x,\eta)\big)  \ot \big(P-\mathscr{V} (x,\eta)\big)\Bigg)
\Bigg[ \mathscr{H}_\eta\left(x,\eta, \frac{1}{2}\big|P-\mathscr{V}(x,\eta)\big|^2 \right)
 \nonumber \\
& -\, \big(P-\mathscr{V} (x,\eta)\big)^\top \mathscr{V}_\eta(x,\eta)  \,
 \mathscr{H}_p \left(x,\eta, \frac{1}{2}\big|P-\mathscr{V}(x,\eta)\big|^2 \right)  \Bigg] . \nonumber
 \end{align}
The naturalness of our structural assumption \eqref{1.13} is also justified by Lagrangian models arising in variational Data assimilation which we describe briefly after the statement of our main result.

\begin{theorem}  \label{theorem1}
Let $\Om\sub \R$ be a bounded open interval and let also
\[
\begin{split}
\mathscr{H}\ :\ \overline{\Om}\by\R^N \by [0,\infty)\larrow [1,\infty),\ \ \ \
\mathscr{V}\ :\ \overline{\Om}\by\R^N \larrow \R^N,
\end{split}
\]
be given maps with $N\in \N$. We suppose that
\beq\label{4.1}
\left\{\ 
\begin{split}
\text{$\mathscr{H}$ is $C^2$ up to the boundary}& , \\
C(|\eta|)\, \geq \, \mathscr{H}_p(x,\eta,p)\ &\geq \ c_0,  \\
2\mathscr{H}_{pp}(x,\eta,p)p\, +\, \mathscr{H}_p(x,\eta,p)\,&\geq\,c_0,   \\ 
 \big| \mathscr{H}_x(x,\eta,p) \big| \,+\,  \big| \mathscr{H}_\eta(x,\eta,p) \big|\, &\leq \, C(|\eta|)(1\, +\, p),      \\
 \big| \mathscr{H}_{pp}(x,\eta,p) \big| \,+\, \big| \mathscr{H}_{p\eta}(x,\eta,p) \big| \,+\, \big| \mathscr{H}_{px}(x,\eta,p) \big| \, &\leq \, C(|\eta|)(1\, +\, p^M), 
\end{split}
\right.
\eeq
and also
\beq \label{4.2}
\left\{\ \ \
\begin{split}
&\text{$\mathscr{V}$ is $C^1$ up to the boundary} , \\
& \big|\mathscr{V}(x,\eta) \big|\, \leq \, (1/c_0) (1\, +\, |\eta|^\al),   
\end{split}
\right.
\eeq
for some constants $c_0,\al \in (0,1)$, some $M\in \N$, some positive continuous increasing function $C \in C^0([0,\infty))$ and all $(x,\eta,p) \in \Om\by \R^N \by [0,\infty)$. Then, for any affine map $b : \R\larrow \R^N$, there exists a map $u^\infty\in W^{1,\infty}_b(\Om,\R^N)$ with the following properties:

\begin{enumerate} 

\item \label{(1)} $u^\infty$ is a vectorial Absolute Minimiser of the functional 
\beq \label{1.17}
E_\infty (u,A)\, =\, \underset{x\in A}{\ess\,\sup}\, \mathscr{H}\left(x,u(x), \frac{1}{2}\big|\D u(x)-\mathscr{V}(x,u(x))
\big|^2 \right),
\eeq
that is it satisfies \eqref{1.1a}.

\item \label{(2)} $u^\infty$ is a $\mD$-solution (see Definitions \ref{definition7} and \ref{definition11}) of the system
\beq \label{1.18}
\ \ \mF_\infty\big(\cdot,u,\D u,\D^2 u\big)\,=\,0,\ \ \ \ \text{on }\Om,
\eeq
where $\mF_\infty$ is given by \eqref{1.14}.

\item \label{(3)} $u^\infty$ is a subsequential limit as $m\ri\infty$ in the strong $W^{1,q}(\Om,\R^N)$ topology of $C^2$ minimisers $\{u^m\}_{m=2}^\infty$ of the functionals
\beq \label{1.19}
 E_m (u,A)\, =\, \int_A \mathscr{H}\left(x,u(x), \frac{1}{2}\big|\D u(x)-\mathscr{V}(x,u(x))
\big|^2 \right)^{\!\! m} dx,
\eeq
where each $u^m$ minimises over the respective space $W^{1,2m}_b(\Om,\R^N)$, for any $q\in [1,\infty)$.

\item \label{(4)} There is an open subset $\Om_\infty\sub \Om$ such that $u^\infty \in C^2(\Om_\infty,\R^N)$. Moreover, 
\[
\ \ \Om\set \Om_\infty\, =\, \p \left(\big\{ \, \D u^\infty = \mathscr{V}(\cdot,u^\infty)  \,\big\}\right) 
\]
and hence $\Om\set \Om_\infty$ is compact and nowhere dense in $\Om$.
\end{enumerate} 

\end{theorem}

Item \eqref{(4)} above is a partial regularity assertion which differs from more classical results in that the singular set $\Om\set \Om_\infty$ is a relatively compact nowhere dense set (a topological boundary) \emph{but not necessarily a Lebesgue nullset}. This is a new type of partial regularity which seems to arise in $L^\infty$. Item \eqref{(3)} indicates the fashion in which $u^\infty$ is obtained, namely via the well-established method of $L^m$ approximations of $L^\infty$ problems as $m\ri \infty$, but also includes a non-trivial fact, the strong convergence of the $L^m$ minimisers $u^m$ together with their first derivatives to $u^\infty$. Note that Theorem \ref{theorem1} above does \emph{not} state that \eqref{(1)} implies \eqref{(2)}, but instead that there is an object $u^\infty$ which satisfies both. {In order to obtain solely \eqref{(1)}, the hypotheses \eqref{4.1}-\eqref{4.2} can be relaxed substantially (accordingly, see the paper \cite{AK}), but since herein we are interested in the satisfaction of the equations as well we do not tackle this problem separately}. Finally, due to the dependence of $\msL$ on the lower order terms $(x,u(x))$, the Absolutely Minimising $\mD$-solution $u^\infty$ is \emph{not} in general unique, as shown by the example $\msL(x,P)=\sin^2x +|P|^2$ of Yu \cite{Y} even when $n=N=1$ and $\mathscr{V}\equiv 0$. Uniqueness is a most delicate question already in the scalar case (see \cite{JWY,MWZ}). Let us also recall that \eqref{(1)} above has been obtained in \cite{BJW1} but under \emph{different hypotheses on $\msL$} which in particular require $\msL(x,\eta,0)=0$ and $\mathscr{V}\equiv 0$, a fact \emph{incompatible} with ``additive" Lagrangians like those arising in Data Assimilation which we describe right next.

{The motivation to study the present $1D$ vectorial $L^\infty$ variational problem comes from Data Assimilation models arising in the Earth Sciences and especially in Meteorology. More precisely, following the terminology of \cite{B}, we are inspired by a continuous time generalisation of what is known as weakly constrained four-dimensional variational assimilation (4D-Var) in the geosciences. For more details we refer e.g.\ to \cite{AJSV, BP, CT, De, FS, LDT, PVT, RCPTV, RJ, TC, Y}.}

{Let us describe briefly the model in pure mathematical terms.} Let $\mathscr{V} : \Om \by \R^N\larrow \R^N$ be a time-dependent vector field describing the law of motion of a body moving along a trajectory defined by the solution  $u:\Om\sub \R \larrow \R^N$ of $\D u=\mathscr{V}(\cdot,u)$ (e.g.\ Newtonian forces, finite-dimensional Galerkin approximation of the Euler equations, etc.). Let also $k: \Om\sub \R \larrow \R^M$ be some partial ``measurements" in continuous time along the trajectory and $K: \R^N \larrow \R^M$ be a submersion which corresponds to some component of the trajectory we are able to measure, for example some projection. Then, we wish to find a $u$ which should satisfy the law of motion and also be compatible with the measurements along the trajectory:
\[
\ \ \ \D u(t) = \mathscr{V}\big(t,u(t)\big) \ \ \ \& \ \ \ K(u(t)) = k(t), \ \ \ t\in \Om.
\]
However, this problem is in general overdetermined (due to errors in the measurements, etc.) since we impose a pointwise constraint to the solution of the system. In standard variational Data Assimilation (see \cite{B, BS}), one instead seeks for approximate solutions by minimising the ``error" integral functional $E_1$ given by \eqref{1.9} for $m=1$ and with Lagrangian given by
\beq \label{1.20}
\mathscr{L}(x,\eta,P)\, :=\, 1+\, \frac{1}{2} \big|k(x)-K(\eta) \big|^2 +\, \frac{1}{2} \big|P-\mathscr{V}(x,\eta) \big|^2
\eeq
which describes the ``error". But if instead we choose to use the respective supremal functional \eqref{1.1} with $\msL$ as in \eqref{1.20}, large ``spikes" of the deviation from the actual solution with small area are from the outset excluded. For the Lagrangian \eqref{1.20}, the equations \eqref{1.5a}-\eqref{1.14} arising in Data Assimilation read
\beq \label{1.21}
\begin{split}
&\big|\D u- \mathscr{V}(\cdot,u)\big|^2  \Big(\D^2 u  - \mathscr{V}_\eta(\cdot,u)-\mathscr{V}_x(\cdot,u) \Big) \, -\, \big[\D u- \mathscr{V}(\cdot,u)\big]^\bot 
\\
& \centerdot  \Big(  \big(K(u)-k\big)^\top K_\eta(u) \, -\, \big(\mathscr{W}u\big)^\top \mathscr{V}_\eta(\cdot,u) \Big)  
\, -\, \Big[  K_\eta(u):\big(K(u)-k\big) 
\\
 & \ot \D u\,  +\,  \big(K(u)-k\big)\cdot k_x   \Big] \big(\D u- \mathscr{V}(\cdot,u)\big)\,=\,0
 \end{split}
\eeq
Our main result applies in particular to \eqref{1.20}-\eqref{1.21}. Although the $L^\infty$ equations are more complicated than the respective $L^2$ Euler-Lagrange equations, evidence obtained in \cite{BK} suggests that they provide more accurate models.

\section{Basics and Young measures into spheres} Our basic notation is either standard  (as e.g.\ in \cite{E2,EG}) or self-explanatory.  For example, the Lebesgue measure on $\R$ will be denoted by $|\cdot|$, the characteristic function of the set $A$ by $\chi_A$, the standard Sobolev and $L^p$ spaces of maps  $u : \Om \sub \R \larrow  \R^N$ by $L^p(\Om,\R^N)$, $W^{m,p}(\Om,\R^N)$, etc. We will also follow the standard practice that while deriving estimates, universal constants may change from line to line but will be denoted by the same letter. $N\in \N$ will always be the dimension of the range of our candidate solutions $u:\Om\sub \R \larrow \R^N$. Unless indicated otherwise, Greek indices $\al,\be,\ga,...$ will run in $\{1,...,N\}$ and the summation convention will be employed in products of repeated indices. The standard basis on $\R^N$ will be denoted by $\{e^1,...,e^N\}$ and hence for the map $u$ with components $u_\al$ we will write $u(x)=u_\al(x)e^\al$.  The norm symbol $|\cdot|$ will always indicate the Euclidean one and the respective inner product will be denoted by ``$\cdot$". Given $\xi \in \R^N$, we define for later use the following orthogonal projections of $\R^N$:
\beq \label{5.4}
[\xi]^\top \,  := \, \sgn(\xi) \ot \sgn(\xi), \ \ \ 
[\xi]^\bot \,  := \, I \, - \, \sgn(\xi) \ot \sgn(\xi). 
\eeq
Here ``sgn" stands for the sign function: $\sgn(\xi):=\xi/|\xi|$ when $\xi \neq 0$ and $\sgn(0):=0$.

Let now $E\sub \R$ be a (Lebesgue) measurable set and consider the Alexandroff $1$-point compactification of the space $\R^N$:
\[
\smash{\overline{\R}}^N\, :=\ \R^N \cup \{\infty\}.
\]
Its topology will be the standard one which makes it isometric to the $N$-sphere (via the stereographic projection which identifies $\{\infty\}$ with the north pole). The space $\R^N$ will be considered equipped with the metric topology induced by the embedding into its compactification $\smash{\overline{\R}}^N$ but balls, distances, etc.\ will be taken with respect to its usual metric structure.

\begin{definition}[Young Measures into the $1$-point compactification of $\R^N$] The space of Young Measures $E\sub \R \larrow \smash{\overline{\R}}^N$ is denoted by $\mY\big(E,\smash{\overline{\R}}^N\big)$ and is the set of probability-valued maps
\[
\R \supseteq E\, \ni \, x\, \lmapsto\, \vartheta(x) \, \in \mathscr{P}\big(\smash{\overline{\R}}^N\big)
\]
which are measurable in the following sense: for any fixed open set $U\sub \smash{\overline{\R}}^N$, the function $E\ni x\lmapsto [\vartheta(x)](U)\in \R$ is measurable. This is called \emph{weak* measurability}.
\end{definition}

The set $\mY\big(E,\smash{\overline{\R}}^N\big)$ is a subset of the unit sphere {of the space $L^\infty_{w^*}\big( E,\mM\big(\smash{\overline{\R}}^N\big) \big)$.} This Banach space consists of weakly* measurable maps valued in the signed Radon measures: $E \ni x \lmapsto \vartheta(x) \in \mM\big(\smash{\overline{\R}}^N\big)$. The norm of the space is
\[
\| \vartheta \|_{L^\infty_{w^*} ( E,\mM(\smash{\overline{\R}}^N) )}\, :=\, \underset{x\in E}{\ess\,\sup}\, \big\|\vartheta(x)\big\|\big(\smash{\overline{\R}}^N\big)
\]
where ``$\|\cdot\|\big(\smash{\overline{\R}}^N\big)$" is the total variation. For more details about this and relevant spaces we refer e.g.\ to \cite{FL} (and references therein). Hence, the Young Measures are the subset of the unit sphere which consists of probability-valued weakly* measurable maps. It can be shown (see e.g.\ \cite{FL}) that $L^\infty_{w^*}\big( E,\mM\big(\smash{\overline{\R}}^N\big) \big)$ is the dual space of the  $L^1$ space of measurable maps valued in the (separable) space $C^0\big(\smash{\overline{\R}}^N\big)$ of real continuous functions over $\smash{\overline{\R}}^N$, in the standard Bochner sense:
\[
\left( L^1\big( E, C^0\big(\smash{\overline{\R}}^N\big)\big) \right)^* \, =\, L^\infty_{w^*}\big( E,\mM\big(\smash{\overline{\R}}^N\big) \big).
\]
The elements of this space are certain Carath\'eodory functions $\Phi : E \by \smash{\overline{\R}}^N \larrow \R$ endowed with the norm
\[
\| \Phi \|_{L^1( E, C^0(\smash{\overline{\R}}^N))}\, :=\, \int_E \big\|\Phi(x,\cdot)\big\|_{C^0(\smash{\overline{\R}}^N)} \, dx.
\]
The space $L^1\big( E, C^0\big(\smash{\overline{\R}}^N\big)\big)$ is separable and the duality pairing 
\[
\langle\cdot,\cdot\rangle\ :\ \  L^\infty_{w^*}\big( E,\mM\big(\smash{\overline{\R}}^N\big) \big) \by L^1\big( E, C^0\big(\smash{\overline{\R}}^N\big)\big) \larrow \R
\]
is given by
\[
\langle \vartheta, \Phi \rangle\, :=\, \int_E \int_{\smash{\overline{\R}}^N} \Phi(x,X)\, d[\vartheta(x)](X)\, dx.
\]
The unit ball of $L^\infty_{w^*}\big( E,\mM\big(\smash{\overline{\R}}^N\big) \big)$ is sequentially weakly* compact. Hence, for any bounded sequence $(\vartheta^m)_1^\infty \sub L^\infty_{w^*}\big( E,\mM\big(\smash{\overline{\R}}^N\big) \big)$, there is a limit map $\vartheta$ and a subsequence of $m$'s along which $\vartheta^m \weakstar \vartheta$ as $m\ri \infty$. 

\begin{remark}[Properties of $\mY(E,\smash{\overline{\R}}^N)$] The set of Young measures is convex and by the compactness of $\smash{\overline{\R}}^N$, it can be proved that it is sequentially weakly* compact in $L^\infty_{w^*}\big( E,\mM\big(\smash{\overline{\R}}^N\big) \big)$ (see e.g.\ \cite{FG, CFV}). This property is essential in our setting. Moreover, the space of  measurable maps $v : E\sub \R \larrow \R^N$ can be (nonlinearly) embedded into $\mY\big(E,\smash{\overline{\R}}^N\big)$ by $v \lmapsto \de_v$, $(\de_v)(x):= \de_{v(x)}$, $x\in E$. 
\end{remark}

The following lemma is a small variant of a standard result about Young measures but it plays an important role in our setting (for the proof see \cite{FG, CFV, V, K8}).

\begin{lemma} \label{lemma2a} Let $v^m,v^\infty : E\sub \R \larrow {\R}^N$ be measurable maps, $m\in \N$. Then, up to the passage to subsequences, the following equivalence holds true as $m\ri \infty$:
\[
\begin{split}
v^m \larrow v^\infty, \text{ a.e.\ on }E \ \ &\Longleftrightarrow \ \ \de_{v^m} \weakstar\, \de_{v^\infty}, \text{ in }\mY\big(E,\smash{\overline{\R}}^N\big).
\end{split}
\]
\end{lemma}

\section{$\mD$-solutions for fully nonlinear systems} 

{
Now we give the main definitions of our notion of solution only in the special case which is needed in this paper. For the general case and applications we refer to \cite{K8}-\cite{K11} and \cite{CKP,KP,KP2}.
\begin{definition}[Diffuse 2nd derivatives]  \label{definition7} {Suppose that $u:\Om\sub \R \larrow \R^N$ is in $W^{1,\infty}_{\text{loc}}(\Om,\R^N)$. For any $h\neq 0$, we consider the difference quotients of the derivative
\[
\D^{1,h}\D u\, =\,  \frac{1}{h}\Big(\D u(\cdot+h)-\D u\Big) \ :\ \, \Om\,\sub \, \R \larrow \R^{N}
\] 
and $\D u$ is understood as being extended by zero on $\R\set \Om$. We define the \textbf{diffuse 2nd derivatives} of $u$  as the subsequential limits $ \mD^2 u $ of $\de_{\D^{1,h}\D u}$ in the space of Young measures from $\Om$ into $\smash{\overline{\R}}^{N}$ along infinitesimal sequences $(h_i)_1^\infty \sub \R\set \{0\}$:}
\[
\begin{split}
\de_{\D^{1,h_{i_j}}\D u} \weakstar \mD^2 u, \ \ \ \ \ \text{ in }\mY\big(\Om,\smash{\overline{\R}}^{N}\big), \ \ \text{ as }i\ri \infty.
\end{split}
\]
\end{definition} 
The weak* compactness of $\mY\big(\Om,\smash{\overline{\R}}^{N}\big)$ implies that every $u\in W^{1,\infty}_{\text{loc}}(\Om,\R^N)$ possesses diffuse $2$nd derivatives, in particular at least one for every choice of $(h_i)_1^\infty$. For our notion of generalised solution, let us first introduce the following notation: if $\vartheta$ is a probability measure on $\smash{\overline{\R}}^{N}$, we define its \emph{reduced support} as
\[
\supp_*(\vartheta)\, :=\, \supp(\vartheta)\set\{\infty\} \ \sub \,\R^N.
\] 
\begin{definition}[$\mD$-solutions of 2nd order ODE systems] \label{definition11} Let 
\beq \label{2.11a}
\mF\ : \ \ (\Om\sub \R) \by \R^N\by \R^{N}\by \R^N \larrow \R^N
\eeq
be a Borel measurable map with $\Om$ an open set. Consider the ODE system
\beq \label{2.11}
\mathcal{F}\big(\cdot,u,\D u,\D^2 u\big)\, =\,0, \ \ \text{ on }\Om.
\eeq
We say that a map $u : \Om\sub \R \larrow \R^N$ in $W^{1,\infty}_{\text{loc}}(\Om,\R^N)$ is a \textbf{$\mD$-solution of \eqref{2.11}} when for any diffuse 2nd derivative $\mD^2u \in \mY\big(\Om,\smash{\overline{\R}}^{N}\big)$ we have
\beq \label{2.12}
\ \ \sup_{X\in \supp_*(\mD^2u(x))}\Big|\mathcal{F}\Big(x,u(x),\D u(x),X\Big)\Big|\, =\,0, \ \ \text{a.e. }x\in\Om.
\eeq 
\end{definition}
In general diffuse derivatives may not be unique for non-differentiable maps. Moreover,  \eqref{2.12} is trivially satisfied at certain points at which it may happen that $\mD^2u(x)=\de_{\{\infty\}}$ and hence $\supp_*(\mD^2u(x)) =\emptyset$ (see also the examples in \cite{K8,K11,K10}). It is an immediate consequence of Lemma \ref{lemma2a} that diffuse derivatives are compatible with classical derivatives, in the sense that if $u$ is twice differentiable a.e.\ on $\Om$, then the diffuse 2nd derivative $\mD^2 u$ is unique in the sense that $\mD^2 u= \de_{\D^2 u}$ a.e.\ on $\Om$. The converse is true as well if $\D^2u$ is interpreted in the sense of Ambrosio-Mal\'y as ``derivative in measure" (see \cite{AM,K8}). As a direct consequence we have that $\mD$-solutions are compatible with a.e.\ twice differentiable strong solutions.}

\section{The fundamental equations arising in $L^\infty$}

In this section we \emph{formally} derive the fundamental equations \eqref{1.5a}-\eqref{1.5c} and in particular \eqref{1.14} associated to $L^\infty$ variational problems for \eqref{1.1}. The formal derivation of \eqref{1.5a}-\eqref{1.5c} has been performed in \cite{K1}, but we include it here because it provides insights of the method of proof which makes the foregoing calculations rigorous. We obtain the $L^\infty$ equations in the limit of the Euler-Lagrange equations related to the $L^m$ integral functional \eqref{1.9} as $m\ri \infty$. Here we suppose that $m\geq 2$. The Euler-Lagrange equation of \eqref{1.9} is the ODE system
\beq \label{5.2}
\D \Big(\mathscr{L}^{m-1}(\cdot , u,\D u)\msL_P(\cdot , u,\D u)\Big) \ = \ \mathscr{L}^{m-1}(\cdot , u,\D u)\msL_\eta (\cdot , u,\D u).
\eeq
By distributing derivatives and normalising, \eqref{5.2} gives
\beq \label{5.3}
\D \big(\msL(\cdot , u,\D u)\big) \msL_P(\cdot , u,\D u) \ + \ \frac{\msL(\cdot , u,\D u)}{m-1}\Big( \D \big(\msL_P(\cdot , u,\D u)\big) - \msL_\eta (\cdot , u,\D u)\Big) \ = \ 0.
\eeq
Then, by employing \eqref{5.4} applied to $\xi=\msL_P(\cdot ,u,\D u)$ we expand the term in the bracket of \eqref{5.3} and obtain
\beq \label{5.5}
\begin{split}    
&\D \big(\msL(\cdot , u,\D u)\big) \msL_P(\cdot , u,\D u) 
\\
& + \ \frac{\msL(\cdot , u,\D u)}{m-1}[\msL_P(\cdot , u,\D u)]^\top 
  \Big(\D \big(\msL_P(\cdot , u,\D u)\big) 
- \msL_\eta (\cdot , u,\D u)\Big)  
\\
=& - \ \frac{\msL(\cdot , u,\D u)}{m-1}[\msL_P(\cdot , u,\D u)]^\bot 
 \Big( \D \big(\msL_P(\cdot , u,\D u)\big) 
- \msL_\eta (\cdot , u,\D u)\Big) .
\end{split}
\eeq
{By mutual orthogonality of the projections in \eqref{5.4}, the left and right hand side of \eqref{5.5} are normal to each other. Hence, they both vanish and we may split the system to two components. We renormalise the second half of \eqref{5.5} by multiplying by $m-1$ to obtain
\beq \label{5.6}
\begin{split}    
\D \big(\msL(\cdot , u,\D u)\big) \msL_P(\cdot , u,\D u) \hspace{155pt}
\\
+ \ \frac{\msL(\cdot , u,\D u)}{m-1}[\msL_P(\cdot , u,\D u)]^\top 
 \Big( \D \big(\msL_P(\cdot , u,\D u)\big) 
- \msL_\eta (\cdot , u,\D u)\Big) \, &=\, 0,
\\
\msL(\cdot , u,\D u) [\msL_P(\cdot , u,\D u)]^\bot \Big( \D \big(\msL_P(\cdot , u,\D u)\big) 
- \msL_\eta (\cdot , u,\D u)\Big)\, &=\,0 .
\end{split}
\eeq
As $m\ri \infty$, we obtain the complete ODE system in $L^\infty$ which after an expansion of derivatives and summation of the systems becomes \eqref{1.5a}-\eqref{1.5c}.}

\smallskip

\noi \textbf{The degenerate elliptic case of the equations in $L^\infty$.} Unfortunately, as we have already explained it is not in general possible to go much further without imposing the structural condition \eqref{1.13} on the hamiltonian $\msL$. The problem is that the system  fails to be degenerate elliptic in the sense needed for existence. In particular, \emph{the coefficient of \eqref{1.5a}-\eqref{1.5c} may be discontinuous at points where $\msL_P(\cdot,u,\D u)=0$}. Assumption \eqref{1.13} forces the matrices $[\msL_P]^\bot$ an $\msL_{PP}$ to commute and also makes the coefficients  \emph{continuous} by allowing them to match after a rescaling.

We now \emph{formally} derive \eqref{1.14} and also the Euler-Lagrange equations of the $L^m$ functional \eqref{1.19} in the expanded form in which it will be used later, under the assumptions \eqref{4.1}-\eqref{4.2}. There is no loss of generality in assuming $\mathscr{H}\geq 1$ since if it is bounded below, we can always add a positive constant to $\mathscr{H}$ and the equations remain the same because additive constants commute with the supremal functional (and this constant also regularises the minimisers of the respective $L^m$ functional). In order to derive the equations, we first differentiate \eqref{1.13} and for the sake of brevity we suppress the argument $(x,\eta,P)$ of $\msL$ and $\big(x,\eta,\frac{1}{2}|P-\mathscr{V}(x,\eta)|^2\big)$ of $\mathscr{H}$ and their respective derivatives:
\begin{align}
\mathscr{L}_{P_\al}\, &=\, \mathscr{H}_p\big(P-\mathscr{V} (x,\eta) \big)_\al , \nonumber
\\
\mathscr{L}_{\eta_\al}\, &=\, \mathscr{H}_{\eta_\al} \, -\, \mathscr{H}_p \big(P-\mathscr{V}(x,\eta)\big)_\ga  \mathscr{V}_{\ga {\eta_\al}}(x,\eta)  , 
\nonumber
\\
\msL_x \, &=\, \, \mathscr{H}_x \, -\, \mathscr{H}_p \big(P-\mathscr{V}(x,\eta)\big)_\ga \mathscr{V}_{\ga  x}(x,\eta) , 
\nonumber
\\
\mathscr{L}_{P_\al P_\be}\!  &=\, \mathscr{H}_{pp} \big(P-\mathscr{V}(x,\eta)\big)_\al \big(P-\mathscr{V}(x,\eta)\big)_\be\, +\, \mathscr{H}_p \,\de_{\al \be},  \nonumber
\\
\mathscr{L}_{P_\al x} \, &= -\mathscr{H}_p \mathscr{V}_{\al x}(x,\eta) \, +\, \big(P-\mathscr{V}(x,\eta)\big)_\al \Big[   \mathscr{H}_{px} \, -\, \mathscr{H}_{pp} \big(P-\mathscr{V} (x,\eta) \big)_\ga \mathscr{V}_{\ga  x}(x,\eta) \Big]   ,
\nonumber
\\
\mathscr{L}_{P_\al \eta_\be} &= -\mathscr{H}_p \mathscr{V}_{\al \eta_\be}(x,\eta)  +  \big(P-\mathscr{V}(x,\eta)\big)_\al  
\Big[   \mathscr{H}_{p \eta_\be} - \mathscr{H}_{pp} \big(P-\mathscr{V}(x,\eta)\big)_\ga \mathscr{V}_{\ga  \eta_\be} (x,\eta) \Big] . 
\nonumber
\end{align}
Recall now that \eqref{5.6} comprises the Euler-Lagrange equations of \eqref{1.9} in expanded form. Since by assumption $\mathscr{H}_p>0$, we have the identities
\beq \label{Id}
\begin{split}
 \big[\mathscr{H}_p\big(P-\mathscr{V}(x,\eta)\big)\big]^\top &=\, \big[P-\mathscr{V}(x,\eta)\big]^\top,
 \\
 \big[\mathscr{H}_p\big(P-\mathscr{V}(x,\eta)\big)\big]^\bot &=\, \big[P-\mathscr{V}(x,\eta)\big]^\bot.
 \end{split}
\eeq
By grouping terms, setting
\beq \label{W(u)}
\mathscr{W}u\,:=\, \D u- \mathscr{V}(\cdot,u)
\eeq 
and omitting the argument $(\cdot,u,\frac{1}{2}|\mathscr{W}u|^2)$ of $\mathscr{H}$ and its derivatives, after a calculation we have
\beq \label{6.4}
\begin{split}
&\Bigg[ \frac{ \mathscr{H} [\mathscr{W}u]^\top}{m-1} \Big(\mathscr{H}_pI \,+\, \mathscr{H}_{pp} (\mathscr{W}u)\ot  (\mathscr{W}u) \Big)   \,  + \,  (\mathscr{H}_p)^2 (\mathscr{W}u)\ot  (\mathscr{W}u)   \Bigg] 
\\
\centerdot &\Big( \D^2 u-\mathscr{V}_\eta(\cdot,u)  \D u -\mathscr{V}_x(\cdot,u) \Big) \, + \, \frac{ \mathscr{H} [\mathscr{W}u]^\top}{m-1} \Bigg(   \mathscr{H}_p(\mathscr{W}u)^\top \mathscr{V}_\eta(\cdot,u) 
 \\
 - & \,  \mathscr{H}_\eta \, +\, \big( \mathscr{H}_{p\eta}\cdot \D u + \mathscr{H}_{px}\big) \,\mathscr{W}u  
 \Bigg) \, +\ \mathscr{H}_p \big(\mathscr{H}_x \,+\, \mathscr{H}_\eta\cdot \D u\big)\, \mathscr{W}u \,=\,0.
\end{split}
\eeq
Similarly, by the identity \eqref{Id} and since $\mathscr{H}>0$, we have
\beq\label{6.5}
\begin{split}
  & [\mathscr{W}u]^\bot \Bigg[  \Big(\mathscr{H}_pI \,+\, \mathscr{H}_{pp} (\mathscr{W}u)\ot  (\mathscr{W}u) \Big) \Big( \D^2 u-\mathscr{V}_\eta(\cdot,u)  \D u -\mathscr{V}_x(\cdot,u) \Big)
    \\
 &- \,  \mathscr{H}_\eta \,+\,  \mathscr{H}_p(\mathscr{W}u)^\top \mathscr{V}_\eta(\cdot,u) \, + \, \big( \mathscr{H}_{p\eta}\cdot \D u\, +\, \mathscr{H}_{px}\big)\, \mathscr{W}u \Bigg]\, =\, 0.
\end{split}
\eeq
Since the projection $[\mathscr{W}u]^\bot $ annihilates $\mathscr{W}u$, \eqref{6.5} simplifies to
\beq\label{6.5A}
\begin{split}
 [\mathscr{W}u]^\bot \left[\mathscr{H}_p \Big( \D^2 u-\mathscr{V}_\eta(\cdot,u)  \D u -\mathscr{V}_x(\cdot,u) \Big) -\mathscr{H}_\eta \, +\, \mathscr{H}_p(\mathscr{W}u)^\top \mathscr{V}_\eta (\cdot,u)  \right]\,=\,0. 
\end{split}
\eeq
We now observe that in view of the identities \eqref{5.4}, the systems \eqref{6.4} and \eqref{6.5A} can be matched and the mutually orthogonal coefficients $[\mathscr{W}u]^\top$ and $[\mathscr{W}u]^\bot$ add to the identity. By multiplying \eqref{6.5A} by $\mathscr{H}_p|\mathscr{W}u|^2$ and adding it to \eqref{6.4}, we obtain
\beq 
\label{6.7}
\begin{split}
  & \Bigg[\frac{ \mathscr{H}\big(\mathscr{H}_p + \mathscr{H}_{pp}|\mathscr{W}u|^2 \big)  [\mathscr{W}u]^\top}{m-1} \, +\,  (\mathscr{H}_p)^2 |\mathscr{W}u|^2I\Bigg] \, \D \big(\mathscr{W}u \big) + \, \frac{ \mathscr{H} [\mathscr{W}u]^\top}{m-1}
\\
  \centerdot & \Bigg( -  \mathscr{H}_\eta \,+ \, \mathscr{H}_p(\mathscr{W}u)^\top \mathscr{V}_\eta(\cdot,u) \, + \, \big( \mathscr{H}_{p\eta}\cdot u + \mathscr{H}_{px}\big)(\mathscr{W}u) \Bigg)\,+\, \mathscr{H}_p \Big(\mathscr{H}_x  
 \\
 & +\, \mathscr{H}_\eta\cdot \D u\Big)\mathscr{W}u\  -\   \mathscr{H}_p |\mathscr{W}u|^2 [\mathscr{W}u]^\bot \Big( \mathscr{H}_\eta \, -\, \mathscr{H}_p(\mathscr{W}u)^\top \mathscr{V}_\eta(\cdot,u)   \Big) \, =\, 0.  
 \end{split}
\eeq
The ODE system \eqref{6.7} is the Euler-Lagrange equation of the functional \eqref{1.19} in expanded form where for sake of brevity we have defined \eqref{W(u)} and suppressed the dependence on the arguments $ (\cdot,u,\frac{1}{2}|\mathscr{W}u|^2 )$ of $\mathscr{H},\mathscr{H}_p,\mathscr{H}_\eta,\mathscr{H}_x$ and $\mathscr{H}_{pp},\mathscr{H}_{p\eta}, \mathscr{H}_{px}$. We also note that the coefficients which are of order $O\big(\frac{1}{m-1}\big)$ are discontinuous, but this causes no problems since the terms involving these will be annihilated as $m\ri \infty$. By letting $m\ri \infty$ in \eqref{6.5A} we obtain \eqref{1.5a} with $\mF_\infty$ given by \eqref{1.14} and $\mathscr{W}u$ by \eqref{W(u)}. We finally rewrite the equations in a form which is more malleable for our proofs later. By setting
\beq
\begin{split} \label{3.21}
\hspace{25pt} F^\infty(\cdot,u,\D u)\, :=& -\, \mathscr{H}_p \Big(\mathscr{H}_x \,+\,\mathscr{H}_\eta\cdot \D u\Big)\mathscr{W}u 
\\
 & \, +\,  (\mathscr{H}_p)^2 |\mathscr{W}u|^2 [\mathscr{W}u]^\bot \Big(  \mathscr{H}_\eta \, -\, \mathscr{H}_p(\mathscr{W}u)^\top \mathscr{V}_\eta(\cdot,u)   \Big),\\
 \end{split} 
 \eeq
 \beq
\label{3.22}
\begin{split}
\ \ \  f^\infty(\cdot,u,\D u)\, :=& - \mathscr{H}[\mathscr{W}u]^\top \Big( -  \mathscr{H}_\eta \,+ \, \mathscr{H}_p(\mathscr{W}u)^\top \mathscr{V}_\eta(\cdot,u) \\
& + \, \big( \mathscr{H}_{p\eta}\cdot u \, + \, \mathscr{H}_{px}\big) \, \mathscr{W}u \Big) ,
\end{split}
\eeq
\beq
 \label{3.22a}
 A^\infty(\cdot,u,\D u)\, := \, \mathscr{H} \big(\mathscr{H}_p \,+\, \mathscr{H}_{pp}|\mathscr{W}u|^2 \big)  [\mathscr{W}u]^\top, \hspace{38pt}
\eeq
\eqref{6.7} can be written as
\beq 
\label{3.23}
\begin{split}
\left[\frac{ A^\infty(\cdot,u,\D u)}{m-1}  +  \mathscr{H}_p^2\Big(\cdot,u,\frac{1}{2}|\mathscr{W}u|^2\Big)  |\mathscr{W}u|^2I \right]  \D \big(\mathscr{W}u \big)   =  F^\infty(\cdot,u,\D u) + \frac{f^\infty(\cdot,u,\D u)}{m-1}
  \end{split}
\eeq
and \eqref{1.5a} as
\beq 
\label{3.24}
\begin{split}
  \mathscr{H}_p^2\Big(\cdot,u,\frac{1}{2}|\mathscr{W}u|^2\Big)  |\mathscr{W}u|^2 \, \D \big(\mathscr{W}u \big)   \,= \,  F^\infty(\cdot,u,\D u)
  \end{split}
\eeq
where $F^\infty,f^\infty, A^\infty$ are given by \eqref{3.21}-\eqref{3.22a} and $\mathscr{W}u$ by \eqref{W(u)}.

\section{Existence of vectorial Absolute Minimisers}

In this section we establish item \eqref{(1)} of our main result Theorem \ref{theorem1} by proving existence of a mapping $u^\infty \in W^{1,\infty}_b(\Om,\R^N)$ which satisfies \eqref{1.1a} for \eqref{1.17}.

\bl[Existence of minimisers and convergence] \label{lemma1} Let $\mathscr{H},\mathscr{V},\Om$ satisfy the assumptions of Theorem \ref{theorem1}. Then,  for any affine mapping $b : \R \larrow \R^N$ and any $m\in \N$, the functional \eqref{1.19} has a minimiser $u^m$ over the space $W^{1,2m}_b(\Om,\R^N)$. Moreover, we have the estimate
\beq \label{6.16}
\| u \|_{W^{1,2m}(\Om)}\ \leq\ C\Big( E_m(u,\Om)^{\frac{1}{2m}}\,+\, \max_{\p\Om}|b|\, +\, 1\Big)
\eeq
for any $u\in W^{1,2m}_b(\Om,\R^N)$, where $C>0$ depends only on the assumptions and the length of $\Om$. In addition, there is a subsequence $(m_k)_1^\infty$ and $u^\infty \in W^{1,\infty}_b(\Om,\R^N)$ such that 
for any $q\geq 1$
\[
\left\{
\begin{split}
 &\ \ \, u^m   \, -\!\!\!\!\larrow  u^\infty,\ \ \  \text{ in }C^0(\overline{\Om},\R^N),\ms\\
&\D u^m  \weak \D u^\infty,\ \text{ in } L^q(\Om,\R^N),
\end{split}
\right.
\]
as $m_k \ri \infty$, and also
\beq
\| u^\infty \|_{W^{1,\infty}(\Om)}\, \leq\, C.
\eeq
Finally, for any $A\sub \Om$ measurable with $|A|>0$, we have the lower semicontinuity inequality
\beq \label{lsc}
E_\infty(u^\infty,A)\, \leq \, \underset{m\ri\infty}{\liminf}\, E_m(u^m,A)^{\frac{1}{m}}.
\eeq
\el

\BPL \ref{lemma1}. \textbf{Step 1.} We begin with some elementary inequalities we use in the sequel. For any $t\geq 0$, $0<\al<1$ and $\e>0$, Young's inequality gives
\beq \label{6.19}
t^\al\, \leq\, \e t\ +\ \Big( \frac{\al}{\e}\Big)^{\frac{\al}{1-\al}}(1-\al).
\eeq
Moreover, for any $P,V \in \R^N$ and $0<\de<1$, we also have
\beq \label{6.20}
(1-\de)|P|^2\, \leq\, |P-V|^2\, +\, \frac{1}{\de}|V|^2.
\eeq
Finally, for any $u \in W^{1,2m}(\Om,\R^N)$, we have the following Poincar\'e inequality \emph{whose constant is uniform in $m\in \N$}:
\beq \label{6.21}
\| u \|_{L^{2m}(\Om)}\, \leq\, 2(|\Om|+1)\Big( \| \D u \|_{L^{2m}(\Om)} \, +\, \max_{\p\Om}|u|\Big).
\eeq
Indeed, in order to see \eqref{6.21}, suppose $u$ is smooth and since $\big| u(x)-u(y) \big| \leq \int_\Om |\D u|$, for $y\in \p\Om$ by H\"older inequality we have
\begin{align}
|u(x)|^{2m}\, &\leq \, \left( \int_\Om |\D u|\, +\, \max_{\p\Om}|u|\right)^{2m} \nonumber\\
&\leq \, 2^{2m-1}\left[ \left( \int_\Om |\D u|\right)^{2m} \, +\, \max_{\p\Om}|u|^{2m} \right]  \nonumber\\
&\leq \, (2(|\Om|+1))^{2m-1} \left[ \int_\Om |\D u|^{2m} \, +\, \max_{\p\Om}|u|^{2m}\right],  \nonumber
\end{align}
which leads to \eqref{6.21}.

\noi \textbf{Step 2.} We now show that the functional $E_m$ is weakly lower semicontinuous in $W^{1,2m}(\Om,\R^N)$. Indeed, by setting
\beq \label{6.22}
H(x,\eta,P)\, :=\, \mathscr{H}\Big(x,\eta, \frac{1}{2}\big|P-\mathscr{V}(x,\eta) \big|^2 \Big)^{\!m},
\eeq
we have for the hessian with respect to $P$ that (we suppress again the arguments of $\mathscr{H}$ and its derivatives)
\[
\begin{split}
H_{PP}(x,\eta,P)\, = &\, m\mathscr{H}^{m-2}\Big[ \mathscr{H} \mathscr{H}_pI\, +\, \Big( \mathscr{H} \mathscr{H}_{pp} \, +\, (m-1)(\mathscr{H}_p)^2\Big) 
\\
&\centerdot\big(P-\mathscr{V}(\cdot,u)\big) \ot \big(P-\mathscr{V}(\cdot,u)\big) \Big].
\end{split}
\]
By \eqref{4.1} and since the projection $[P-\mathscr{V}(\cdot,u)]^\top$ satisfies the matrix inequality $[P-\mathscr{V}(\cdot,u)]^\top\! \leq I$, we obtain 
\begin{align}
H_{PP}(x,\eta,P)\, &\geq\, m\mathscr{H}^{m-2}\Big[ \mathscr{H} \mathscr{H}_pI\, +\, \mathscr{H} \mathscr{H}_{pp}\big(P-\mathscr{V}(\cdot,u)\big) \ot \big(P-\mathscr{V}(\cdot,u)\big) \Big]  \nonumber\\
&\geq \, m\mathscr{H}^{m-1}\Big(\mathscr{H}_p I\, +\, (c_0-\mathscr{H}_p) \big[P-\mathscr{V}(\cdot,u)\big]^{\!\top}\Big)  
\\
&\geq \, m \Big(\mathscr{H}_p I\, +\, (c_0-\mathscr{H}_p) \big[P-\mathscr{V}(\cdot,u)\big]^{\!\top}\Big) \nonumber \\
& \geq  \, mc_0I.  \nonumber
\end{align}
The conclusion now follows by standard lower semicontinuity results (e.g.\ \cite{D, GM}).

\noi \textbf{Step 3.} Now we derive the energy estimate which guarantees the coercivity of $E_m$. By our assumptions on $\mathscr{H}$, there is a $\hat{p}\in [0,p]$ such that
\[
\mathscr{H}(x,\eta,p)\, =\, \mathscr{H}_p(x,\eta,\hat{p})p\, +\, \mathscr{H}(x,\eta,0)\, \geq \, c_0 p\, +\,1.
\]
Hence, by using \eqref{6.20} the above gives
\beq \label{6.23}
\mathscr{H}\Big(x,\eta, \frac{1}{2}\big|P-\mathscr{V}(x,\eta) \big|^2 \Big)\, \geq \, \frac{c_0}{2}(1-\de)|P|^2\, -\, \frac{c_0}{2\de}|\mathscr{V}(x,\eta)|^2. 
\eeq
Then, by \eqref{4.2} and \eqref{6.19}-\eqref{6.20}, for $\si>0$ small we have
\begin{align} \label{6.24}
\mathscr{H}\Big(x,\eta, \frac{1}{2}\big|P-\mathscr{V}(x,\eta) \big|^2 \Big)\, & \geq \, \frac{c_0}{2}(1-\de)|P|^2\, -\, \frac{1}{2c_0\de} (1\, +\, |\eta|^{\al})^2 \nonumber\\
& \geq \, \frac{c_0}{2}(1-\de)|P|^2\, -\, \frac{\si}{c_0\de} |\eta|^2 \, -\, C(\de,\si,\al) ,\nonumber
\end{align}
where $ C(\de,\si,\al)$ denotes a constant depending only on the numbers $\de,\si,\al$. We now select $\de := 1/2$, $\si := 2c_0\e >0$ to find
\[
\mathscr{H}\Big(x,\eta, \frac{1}{2}\big|P-\mathscr{V}(x,\eta) \big|^2 \Big)\, \geq \, \frac{c_0}{4}|P|^2\, -\, \e |\eta|^2 \, -\, C(\e,\al).
\]
Hence, for any $m\in \N$ by the H\"older inequality and the above estimate, we have
\[
\frac{1}{3^{m-1}}\Big(\frac{c_0}{4}\Big)^m |P|^{2m}\, \leq\ \mathscr{H}\Big(x,\eta, \frac{1}{2}\big|P-\mathscr{V}(x,\eta) \big|^2 \Big)^{\! m}\, +\, \e^m |\eta|^{2m} \, +\, C(\e,\al)^{2m}.
\]
Consequently, for any $u\in W^{1,2m}_b(\Om,\R^N)$, by integrating over $\Om$ and by utilising \eqref{6.21} and \eqref{1.19}, we deduce
\[
\begin{split}
3\Big(\frac{c_0}{12}\Big)^m \int_\Om |\D u|^{2m} \, &\leq \, E_m(u,\Om)\, +\, \e^m \int_\Om |u|^{2m} \, +\, C(\e,\al)^{2m}|\Om| \nonumber\\
&\leq \, E_m(u,\Om)\, +\, C(\e,\al)^{2m}|\Om| \\
& \ \ \ \ +\, \e^m  \big(2(|\Om|+1)\big)^{2m} \left\{ \max_{\p\Om}|b|^{2m} \, +\, \int_\Om |\D u|^{2m} \right\}. \nonumber 
\end{split}
\]
Hence, we have obtained the estimate
\[
\left\{ \Big(\frac{c_0}{12}\Big)^m -\, \big(4(|\Om|+1)^2\e\big)^m  \right\} \int_\Om |\D u|^{2m} \, \leq \, E_m(u,\Om)\,  +\, C^{2m}\Big( \max_{\p\Om}|b|^{2m} \, +\, 1\Big)
\]
where $C$ above depends on $\e,\al,\Om$. By choosing 
$\e := c_0/\big(3 \, 2^{5} (|\Om|+1)^{2})$, we get
\[
\left\{ \frac{c_0}{12}  \Big(1-\frac{1}{2^m} \Big)^{\frac{1}{m}}\right\}^m 
\int_\Om |\D u|^{2m} \, \leq \, E_m(u,\Om)\,  +\, C^{2m}\Big( \max_{\p\Om}|b|^{2m} \, +\, 1\Big)
\]
and since $\lim_{m\ri \infty}\big(1- 2^{-m} \big)^{1/m}=1$, the desired estimate \eqref{6.16} ensues.

\noi \textbf{Step 4.} {We show existence of minimisers and convergence by using ideas of \cite{BJW1}}. We have the a priori energy bounds (recall the notation \eqref{W(u)})
\begin{align} 
\inf\Big\{ E_m(v,\Om)^{\frac{1}{2m}} \, :\, v \in &W^{1,2m}_b(\Om,\R^N)\Big\} \, \leq \, E_m(b,\Om)^{\frac{1}{2m}}  \nonumber\\
&=\, \left(  \int_\Om \mathscr{H}\Big(\cdot,b, |\mathscr{W}b|^2/2 \Big)^{\! m}\right)^{\frac{1}{2m}}  \nonumber\\
&\leq\, |\Om|^{\frac{1}{2m}}
\Big\|  \mathscr{H}\Big(\cdot,b, |\mathscr{W}b|^2 / 2\Big)\Big\|^{\frac{1}{2}}_{L^\infty(\Om)}  \nonumber
\end{align}
and $E_m(v,\Om)\geq  0$, for any $v \in W^{1,2m}_b(\Om,\R^N)$. Hence, by standard results (see e.g.\ \cite{D, GM}), there exists a minimiser $u^m$ of the functional $E_m$ in $W^{1,2m}_b(\Om,\R^N)$. Moreover, by \eqref{6.16} and \eqref{6.26} we have the bound
\beq \label{6.25}
\|u^m\|_{W^{1,2m}(\Om)} \, \leq\, C\left( \sup_{\Om} \mathscr{H}\Big(\cdot,b,  |\mathscr{W}b|^2/2 \Big) ^{\!\frac{1}{2}}\,  +\, \max_{\p\Om}|b|\, +\, 1 \right).
\eeq
Let $C(\Om,b)$ denote the right hand side of \eqref{6.25}. Then, for any $q\in [2,m]$, we have
\beq
\begin{split} \label{6.26}
\|u^m\|_{W^{1,2q}(\Om)} \, &\leq\, |\Om|^{\frac{1}{2q}-\frac{1}{2m}} \|u^m\|_{W^{1,2m}(\Om)} \\
& \leq\, |\Om|^{\frac{1}{2r}-\frac{1}{2m}} C(\Om,b).  
\end{split}
\eeq
Hence, for any $q\geq 1$ fixed, the sequence $(u^m)_1^\infty$ is bounded in $W^{1,2q}_b(\Om,\R^N)$. As such, there exists $u^\infty \in \cap_{q=1}^\infty W^{1,2q}_b(\Om,\R^N)$ satisfying $u^m \weak u^\infty$ in $W^{1,2q}_b(\Om,\R^N)$ along a subsequence $m_k\ri \infty$. By letting $m\ri \infty$ in \eqref{6.26} along the subsequence, the weak lower semicontinuity of the $L^{2q}(\Om,\R^N)$ norm implies 
\[
\|u^\infty\|_{W^{1,2q}(\Om)} \,\leq \, |\Om|^{\frac{1}{2q}} C(\Om,b). 
\]
By letting now $q\ri \infty$, we derive the desired bound for $u^\infty$.

\noi \textbf{Step 5}.  {We finally show \eqref{lsc} by using ideas of \cite{BJW1}}. Fix $A\sub \Om$ a measurable set with $|A|>0$. By recalling that $u^m \weak u^\infty$ as $m\ri \infty$ along a subsequence in $L^q(A,\R^N)$ for any $q\geq 1$, by weak lower semicontinuity we have
\[
\begin{split}
E_\infty(u^\infty,A)\, &=\, \underset{q\ri\infty}{\lim}\, E_q(u^\infty,A)^{\frac{1}{q}} 
\\
& \leq \, \underset{q\ri\infty}{\liminf}\, \left(\underset{m\ri\infty}{\liminf}\, E_q(u^m,A)^{\frac{1}{q}}\right) 
\\
& \leq \, \underset{q\ri\infty}{\liminf}\, \left(\underset{m\ri\infty}{\liminf}\, |A|^{\frac{1}{q}-\frac{1}{m}}\, E_m(u^m,A)^{\frac{1}{m}} \right) 
\\
& \leq\, \underset{m\ri\infty}{\liminf}\, E_m(u^m,A)^{\frac{1}{m}} .
\end{split}
\]
The lemma ensues.       \qed

\ms

Now we are ready to establish the existence of Absolute Minimisers of the functional \eqref{1.17}.

\ms

\noi \textbf{Proof of \eqref{(1)} of Theorem \ref{theorem1}}. We continue from the proof of the previous lemma. Fix $\Om' \Subset \Om$. Since $\Om'$ is a countable disjoint union of intervals, it suffices to assume that $\Om'=(a,b) \Subset \Om$. We fix $\phi \in W^{1,\infty}_0\big((a,b),\R^N\big)$ and set $\psi^\infty := u^\infty +\phi$. Hence, in order to show that $u^\infty$ is an Absolute Minimiser of \eqref{1.17} over $\Om$ it suffices to show that
\[
E_\infty \big(u^\infty,(a,b)\big) \, \leq\, E_\infty \big(\psi^\infty,(a,b)\big).
\]
Note also that $u^\infty(a)=\psi^\infty(a)$ and $u^\infty(b)=\psi^\infty(b)$. We now fix $0<\ga,\de<(b-a)/{3}$ and define the following map:
\[
\psi^{m,\ga,\de}(x)\,:=\, 
\left\{
\begin{array}{ll}
\left(\dfrac{(a+\ga)-x}{\ga}\right) u^m(a) + \left(\dfrac{x-a}{\ga}\right)\psi^\infty(a+\ga), &   x\in (a,a+\ga), \ms\\
\ \psi^\infty(x), & x\in (a+\ga,b-\de), \ms\\
\left(\dfrac{b-x}{\de} \right)\psi^\infty(b-\de) + \left(\dfrac{x-(b-\de)}{\de}\right)u^m(b), &   x\in (b-\de,b) ,
\end{array}
\right.
\]
where $m\in \N \cup\{\infty\}$. Then, we have $\psi^{m,\ga,\de} \in W^{1,\infty}_{u^m}\big((a,b),\R^N\big)$ and 
\[
\D \psi^{m,\ga,\de}(x)\,=\, 
\left\{
\begin{array}{ll}
 \dfrac{\psi^\infty(a+\ga) - u^m(a)}{\ga}, &  \text{on } (a,a+\ga) \ms\\
 \D \psi^\infty, &  \text{on } (a+\ga,b-\de) \ms \ms\\
 \dfrac{ \psi^\infty(b-\de) - u^m(b) }{-\de}, &   \text{on } (b-\de,b) 
\end{array}
\right.
\]
Note now that 
\beq \label{5.11}
\text{$\psi^{m,\ga,\de} \larrow \psi^{\infty,\ga,\de}$ in $W^{1,\infty}\big((a,b),\R^N\big)$, as $m\ri \infty$.}
\eeq
Indeed, since obviously $\psi^{m,\ga,\de} \larrow \psi^{\infty,\ga,\de}$ in $L^\infty\big((a,b),\R^N\big)$, it suffices to note that for a.e.\ $x\in (a,b)$
\[
\begin{split}
\Big| \D\psi^{m,\ga,\de}(x) - \D\psi^{\infty,\ga,\de}(x)\Big|\, &=\, \chi_{(a,a+\ga)}\frac{|u^\infty(a)-u^m(a)|}{\ga} \,+\, \chi_{(b-\de,b)}\frac{|u^\infty(b)-u^m(b)|}{\de}
\\
& \leq\, \left(\frac{1}{\ga} +\frac{1}{\de}\right) \|u^m-u^\infty\|_{L^\infty(\Om)}
\end{split}
\]
and hence
\[
\begin{split}
\Big\| \D\psi^{m,\ga,\de} - \D\psi^{\infty,\ga,\de}\Big\|_{L^\infty(\Om)} \leq\, \left(\frac{1}{\ga} +\frac{1}{\de}\right) \|u^m-u^\infty\|_{L^\infty(\Om)}\, \larrow 0,
\end{split}
\]
as $m\ri \infty$ along a subsequence. Since for each $m\in \N$ $u^m$ is a minimiser of \eqref{1.19} over $W^{1,2m}_b(\Om,\R^N)$, by recalling that $\psi^{m,\ga,\de}=u^m$ at $\{a,b\}$, minimality and H\"older inequality give
\beq \label{5.12}
\begin{split}
E_m \big(u^m,(a,b)\big)^{\frac{1}{m}} \, &\leq \, E_m \big(\psi^{m,\ga,\de},(a,b)\big)^{\frac{1}{m}}
\\
&\leq \, (b-a)^{\frac{1}{m}}E_\infty \big(\psi^{m,\ga,\de},(a,b)\big) .
\end{split} 
\eeq
On the other hand, since 
\[
\begin{split}
E_\infty \big(\psi^{m,\ga,\de},(a,b)\big)\, =& \, \max\Big\{ E_\infty  \big(\psi^{m,\ga,\de},
(a,a+\ga)\big), \\
&\ \ \ \ \ \  \ \ \ E_\infty  \big(\psi^{m,\ga,\de},(a+\ga,b-\de)\big),
\\
&\ \ \ \ \ \  \ \ \ E_\infty  \big(\psi^{m,\ga,\de},(b-\de,b)\big)\Big\} 
\end{split}
\]
and $\psi^{m,\ga,\de}=\psi^\infty$ on $(a+\ga,b-\de)$, we have
\beq  \label{5.13}
\begin{split}
E_\infty \big(\psi^{m,\ga,\de},(a,b)\big)\, 
 \leq & \, \max\Big\{ E_\infty  \big(\psi^{m,\ga,\de},
(a,a+\ga)\big), \,  E_\infty  \big(\psi^\infty,(a ,b )\big),
\\
& \ \ \ \ \ \  \ \ \ E_\infty  \big(\psi^{m,\ga,\de},(b-\de,b)\big)\Big\}.
\end{split}
\eeq
{By combining \eqref{5.11}-\eqref{5.13} and \eqref{lsc} for $A=(a,b)$, we get along a subsequence $(m_i)_1^\infty$ that
\[
\begin{split}
E_\infty \big(u^\infty,(a,b)\big)\, \leq\, \liminf_{i\ri\infty}\Big( \max\Big\{ & E_\infty  \big(\psi^{m_i,\ga,\de},
(a,a+\ga)\big), \,  E_\infty  \big(\psi^\infty,(a ,b )\big),
\\
& E_\infty  \big(\psi^{m_i,\ga,\de},(b-\de,b)\big)\Big\} \Big)
\end{split}
\]
which in turn gives}
\beq \label{5.14}
\begin{split}
E_\infty \big(u^\infty,(a,b)\big)\,   \leq\, \max\Big\{ & E_\infty  \big(\psi^\infty,(a ,b )\big),\,  E_\infty  \big( \psi^{\infty, \ga,\de}  ,
(a,a+\ga)\big) ,
\\
&  E_\infty  \big(  \psi^{\infty,\ga,\de},(b-\de,b)\big)\Big\} \Big).
\end{split}
\eeq
Further, since 
\[
\begin{split}
&\text{$\D\psi^{\infty,\ga,\de} \,\equiv \,\D^{1,\ga}\psi^\infty(a)$, \ \ \ on $(a,a+\ga)$}, 
\\
&\text{$\D\psi^{\infty,\ga,\de}\,\equiv \,\D^{1,-\de}\psi^\infty(b)$, \ \  on $(b-\de,b)$,
}
\end{split}
\]
we have
\beq \label{5.15}
\left\{\ \ \
\begin{split}
E_\infty  \big(\psi^{\infty,\ga,\de},
(a,a+\ga)\big) \, &=\, \max_{[a,a+\ga]} \msL\Big(\cdot, \psi^{\infty,\ga,\de},\D^{1,\ga}\psi^\infty(a)\Big),
\\
E_\infty  \big(\psi^{\infty,\ga,\de},
(b-\de,b)\big) \, &=\, \max_{[b-\de,b]} \msL\Big(\cdot, \psi^{\infty,\ga,\de},\D^{1,-\de}\psi^\infty(b)\Big).
\end{split}
\right.
\eeq
By \eqref{5.14}-\eqref{5.15}, we see that it suffices to show that there exist infinitesimal sequences $(\ga_i)_{i=1}^\infty$ and $(\de_i)_{i=1}^\infty$ such that
\beq \label{5.16}
\begin{split}
E_\infty  \big(\psi^\infty ,
(a,b)\big) \, \geq \, \max\Big\{ \limsup_{i\ri\infty}  &
\max_{[a,a+\ga_i]} \msL\Big(\cdot, \psi^{\infty,\ga_i,\de_i},\D^{1,\ga_i}\psi^\infty(a)\Big),
\\
\limsup_{i\ri\infty}  & \max_{[b-\de_i,b]} \msL\Big(\cdot, \psi^{\infty,\ga_i,\de_i},\D^{1,-\de_i}\psi^\infty(b)\Big) \Big\}.
\end{split}
\eeq
The rest of the proof is devoted to establishing \eqref{5.16} and this will complete the proof. 

We illustrate the idea of the proof of \eqref{5.16} by \emph{assuming first that $\D \psi^\infty$ exists on $[a,b]$}. In this special case, we have
\[
\begin{split}
E_\infty  \big(\psi^\infty ,
(a,b)\big) \, = \,  \sup_{[a,b]} \, & \msL \Big(\cdot, \psi^{\infty},\D\psi^\infty\Big)
\\
\geq \, \max & \Big\{ \msL\Big(a, \psi^{\infty}(a),\D\psi^\infty(a)\Big), 
\,  \msL\Big(b, \psi^{\infty}(b),\D\psi^\infty(b)\Big) \Big\}.
\end{split}
\]
Further, since the difference quotients satisfy
\[
\D^{1,\ga}\psi^\infty(a) \larrow \D\psi^\infty(a), \ \ \ \D^{1,-\de}\psi^\infty(b) \larrow \D\psi^\infty(b),
\]
as $\ga,\de\ri0$, and also
\beq \label{5.17}
\left\{\ \ \ 
\begin{split}
& \max_{a\leq x\leq a+\ga} \Big| \psi^{\infty,\ga,\de}(x) -\psi^\infty(a)\Big| \, \larrow 0, \ \ \text{ as }\ga \ri 0,
\\
& \ \max_{b-\de\leq x\leq b} \Big| \psi^{\infty,\ga,\de}(x) -\psi^\infty(b)\Big| \, \larrow 0, \ \ \text{ as }\de \ri 0,
\end{split}
\right.
\eeq
we obtain
\[
\lim_{\ga\ri 0}  
\max_{[a,a+\ga]} \msL\Big(\cdot, \psi^{\infty,\ga,\de},\D^{1,\ga}\psi^\infty(a)\Big) \, =\, \msL\Big(a, \psi^{\infty}(a),\D\psi^\infty(a)\Big)  ,
\]
\[
\lim_{\de\ri0}   \max_{[b-\de,b]} \msL\Big(\cdot, \psi^{\infty,\ga,\de},\D^{1,-\de}\psi^\infty(b)\Big) \, =\, \msL\Big(b, \psi^{\infty}(b),\D\psi^\infty(b)\Big).
\]
By putting these together we are led to \eqref{5.16}.

Now we return to the general case. Fix $u\in W^{1,\infty}(\Om,\R^N)$, $x\in [a,b]$ and $\e>0$ small and set 
\[
A_\e(x)\,:=\, [x-\e,x+\e] \cap [a,b].
\]
Then, we claim that there is an increasing modulus of continuity $\om \in C^0(0,\infty)$ with $\om(0^+)=0$ such that
\beq \label{5.18}
E\big(u,A_\e(x) \big)\, \geq\, \underset{y \in A_\e(x)}{\ess\,\sup}\ \msL\Big(x,u(x),\D u(y) \Big)\, -\, \om(\e).
\eeq 
In order to see \eqref{5.18}, note that for a.e.\ $y \in A_\e(x)$ we have $|x-y|\leq \e$ and hence by the continuity of $\msL$ and the essential boundedness of $Du$, there is an $\om$ such that
\[
\Big| \msL\Big(x,u(x),\D u(y) \Big) - \msL\Big(y,u(y),\D u(y) \Big) \Big| \, \leq\, \om(\e),
\]
for a.e.\ $y \in A_\e(x)$. Hence,  \eqref{5.18} ensues. Now we claim that
\beq  \label{5.19}
\sup_{A_\e(x)} \left(\limsup_{t\ri0} \frac{1}{2}\Big| \D^{1,t}u-\mathscr{V}(\cdot,u)\Big|^2\right)\, \leq\, 
\underset{A_\e(x)}{\ess\,\sup} \, \frac{1}{2}\Big| \D u-\mathscr{V}(\cdot,u)\Big|^2.
\eeq
In order to see \eqref{5.19}, it suffices to apply the inequality
\[
\left|\frac{v(y+t)-v(y)}{t} \right|\, \leq\, \underset{A_\e(x)}{\ess\,\sup} \, |\D v|, \quad  y,y+t\in A_\e(x),\ t\neq0,
\] 
to the Lipschitz map 
\[
v(y)\,:=\, u(y)\,-\, \int_a^y\mathscr{V}(t,u(t))\,dt
\]
and note the identities
\[
\begin{split}
\D v(y) \, & =\, \D u(y)\,-\, \mathscr{V}(y,u(y)), \ \ \ \ \text{a.e. }y \in A_\e(x),
\\
\ \ \ \ \ \ \ \ \ \D^{1,t}v(y)\,&=\, \D^{1,t}u(y)\,-\, \frac{1}{t}\int_y^{y+t}\mathscr{V}(t,u(t))\,dt
\\
&=\, \D^{1,t}u(y)\,-\, \mathscr{V}(y,u(y)) \, +\,o(1),\ \ \ \  \text{ as }t\ri0,\ y \in A_\e(x).
\end{split}
\]
Hence, \eqref{5.19} holds true. Now, we combine \eqref{1.13}, \eqref{5.18} and \eqref{5.19} together with \eqref{4.1} to find that for any fixed $u\in W^{1,\infty}(\Om,\R^N)$, $x\in [a,b]$ and $\e>0$ small we have
\[
\begin{split}
E_\infty\big(u,(a,b)\big)\, &\geq\, E_\infty\big(u,A_\e(x)\big) \\
&=\, \underset{A_\e(x)}{\ess\,\sup}\ \mathscr{H}\left(y,u(y),\frac{1}{2}\big|\D u-\mathscr{V}(y,u(y))\big|^2 \right)
\\
&\geq \, \underset{y\in A_\e(x)}{\ess\,\sup}\ \mathscr{H}\left(x,u(x),\frac{1}{2}\big|\D u-\mathscr{V}(y,u(y))\big|^2 \right) \, -\, \om(\e)
\\
&= \, \mathscr{H}\left(x,u(x), \underset{y\in A_\e(x)}{\ess\,\sup}\, \frac{1}{2}\big|\D u(y)-\mathscr{V}(y,u(y))\big|^2 \right) \, -\, \om(\e)
\\
&\geq \, \mathscr{H}\left(x,u(x), \underset{y\in A_\e(x)}{\sup}\left[\limsup_{t\ri0}\, \frac{1}{2}\big|\D^{1,t} u(y)-\mathscr{V}(y,u(y))\big|^2\right] \right) \, -\, \om(\e)
\\
&= \, \underset{y\in A_\e(x)}{\sup}\left[\limsup_{t\ri0}\, \mathscr{H}\left(x,u(x), \frac{1}{2}\big|\D^{1,t} u(y)-\mathscr{V}(y,u(y))\big|^2 \right)\right] \, -\, \om(\e)
\\
&\geq \ \limsup_{t\ri0}\, \mathscr{H}\left(x,u(x), \frac{1}{2}\big|\D^{1,t} u(x)-\mathscr{V}(x,u(x))\big|^2 \right) \, -\, \om(\e)
\end{split}
\]
and by letting $\e\ri 0$, we get
\beq \label{5.20}
\ \ E_\infty\big(u,(a,b)\big)\, \geq \ \limsup_{t\ri0}\, \mathscr{H}\left(x,u(x), \frac{1}{2}\big|\D^{1,t} u(x)-\mathscr{V}(x,u(x))\big|^2 \right),
\eeq
for any fixed $u\in W^{1,\infty}(\Om,\R^N)$ and $x\in [a,b]$. Note now that since 
\[
\big|\D^{1,t}u(x)\big| \, \leq \, \|\D u\|_{L^\infty(\Om)}, \ \ \ \ x\in (a,b),\ t\neq 0,
\]
for any infinitesimal sequence $(t_i(x))_{i=1}^\infty$ there is a subsequence denoted again by the same symbol such that 
\beq \label{5.21}
\text{the limit }\ \lim_{i\ri \infty}\D^{1,t_i(x)}u(x)\  \text{ exists in }\R^N.
\eeq
By \eqref{5.20}-\eqref{5.21} and the continuity of $\mathscr{H}$ we find that
\beq \label{5.22}
\begin{split}
E_\infty\big(u,(a,b)\big)\, &  \geq \ \limsup_{i\ri \infty}\, \mathscr{H}\left(x,u(x), \frac{1}{2}\big|\D^{1,t_i(x)} u(x)-\mathscr{V}(x,u(x))\big|^2 \right) 
\\
&= \,  \mathscr{H}\left(x,u(x), \frac{1}{2}\Big| \lim_{i\ri \infty}\D^{1,t_i(x)}u(x)-\mathscr{V}(x,u(x))\Big|^2 \right).
\end{split}
\eeq
Now we apply \eqref{5.22} to
\[
u\,=\, \psi^\infty,\ \ x\,=\, a,\, b
\]
to infer that there exist sequences $(\ga_i)_{i=1}^\infty$ and $(\de_i)_{i=1}^\infty$ such that 
\beq \label{5.23a}
\text{ the limits }\ \lim_{i\ri \infty}\D^{1,\ga_i}\psi^\infty(a) \ \ \text{ and } \ \lim_{i\ri \infty}\D^{1,-\de_i}\psi^\infty(b) \ \ \text{ exist in $\R^N$}
\eeq
and also
\begin{align} 
\label{5.23}
E_\infty\big(\psi^\infty,(a,b)\big)\, \geq \, \max\Bigg\{ & \mathscr{H}\left(a,\psi^\infty(a), \frac{1}{2}\Big|\lim_{i\ri \infty}\D^{1,\ga_i}\psi^\infty(a)-\mathscr{V}(a,\psi^\infty(a))\Big|^2 \right)
,\nonumber 
\\
& \mathscr{H}\left(b,\psi^\infty(b), \frac{1}{2}\Big|\lim_{i\ri \infty}\D^{1,-\de_i}\psi^\infty(b)-\mathscr{V}(b,\psi^\infty(b))\Big|^2 \right)
\Bigg\}.
\end{align}
On the other hand, by \eqref{5.15}, \eqref{5.17} and \eqref{5.23a}, for $\ga=\ga_i$ and $\de=\de_i$ we have
\begin{align}
\label{5.24}
\lim_{i\ri\infty} E_\infty  \big(\psi^{\infty,\ga_i,\de_i},
(a,a+\ga_i)\big) \, &=\, \lim_{i\ri\infty} \max_{[a,a+\ga_i]} \msL\Big(\cdot, \psi^{\infty,\ga_i,\de_i},\D^{1,\ga_i}\psi^\infty(a)\Big)
\nonumber 
\\
=&\, \mathscr{H}\left(a,\psi^\infty(a), \frac{1}{2}\Big|\lim_{i\ri \infty}\D^{1,\ga_i}\psi^\infty(a)-\mathscr{V}(a,\psi^\infty(a))\Big|^2 \right)
\end{align}
and similarly
\begin{align}
\label{5.25}
\lim_{i\ri\infty} E_\infty  \big(\psi^{\infty,\ga_i,\de_i},
(b-\de_i,b)\big) \, &=\, \lim_{i\ri\infty} \max_{[b-\de_i,b]} \msL\Big(\cdot, \psi^{\infty,\ga_i,\de_i},\D^{1,-\de_i}\psi^\infty(b)\Big)
\nonumber 
\\
=&\, \mathscr{H}\left(b,\psi^\infty(b), \frac{1}{2}\Big|\lim_{i\ri \infty}\D^{1,-\de_i}\psi^\infty(b)-\mathscr{V}(b,\psi^\infty(b))\Big|^2 \right)
\end{align}
By putting together \eqref{5.23}-\eqref{5.25} we see that \eqref{5.16} ensues and so does item \eqref{(1)} of Theorem \ref{theorem1}.    \qed \ms

\ms

\section{Existence of $\mD$-solutions to the equations in $L^\infty$}

In this section we establish items \eqref{(2)}-\eqref{(4)} of Theorem \ref{theorem1}. 
We begin by showing that the minimisers obtained in the previous section actually are weak solutions of the respective Euler-Lagrange equations.

\bl[Weak solutions of the $L^m$ equations] \label{lemma2} Let $\mathscr{H},\mathscr{V},\Om,b$ satisfy the assumptions of Theorem \ref{theorem1} and let $(u^m)_1^\infty$ be the sequence of minimisers constructed in Lemma \ref{lemma1}. Then, each $u^m$ is a weak solution in  $W^{1,2m}_b(\Om,\R^N)$ of the Euler-Lagrange equations of \eqref{1.19} on $\Om$:
\beq \label{6.27}
\begin{split}
\ \ \ \ \ \ & \D \left(\mathscr{H}^{m-1}\Big(\cdot,u,\dfrac{1}{2}\big| \mathscr{W}u\big|^2\Big)\mathscr{H}_p\Big(\cdot,u,\dfrac{1}{2}\big| \mathscr{W}u\big|^2\Big)\, \mathscr{W}u\right) \\
  & =\ \mathscr{H}^{m-1}\Big(\cdot,u,\dfrac{1}{2}\big| \mathscr{W}u\big|^2\Big)\mathscr{H}_p\Big(\cdot,u,\dfrac{1}{2}\big| \mathscr{W}u\big|^2\Big)
\\
& \ \ \ \centerdot \left( \mathscr{H}_\eta\Big(\cdot,u,\dfrac{1}{2}\big| \mathscr{W}u\big|^2\Big)\, -\, \mathscr{H}_p\Big(\cdot,u,\dfrac{1}{2}\big| \mathscr{W}u\big|^2\Big)(\mathscr{W}u )^\top \mathscr{V}_\eta(\cdot,u)\right),
\end{split}
\eeq
where $\mathscr{W}u$ is given by \eqref{W(u)}.
\el

\BPL \ref{lemma2}. {Let $H$ be given by \eqref{6.22}. By suppressing for brevity the arguments of $\mathscr{H},\mathscr{H}_p,\mathscr{H}_\eta$, we have 
\[
\begin{split}
H_P(x,\eta,P)\, &=\, m\mathscr{H}^{m-1}\mathscr{H}_p\, \big(P-\mathscr{V}(x,\eta)\big),\\
H_\eta(x,\eta,P) \, &=\, m\mathscr{H}^{m-1}\mathscr{H}_p\, \Big( \mathscr{H}_\eta\,-\, \mathscr{H}_p \big(P-\mathscr{V}(x,\eta)\big)^\top \mathscr{V}_\eta(x,\eta)\Big)
\end{split}
\]
and the ODE system \eqref{6.27} can be written compactly as 
\beq \label{compressed}
\D \big(H_P(\cdot,u,\D u)\big)\,=\, H_\eta(\cdot,u,\D u).
\eeq
By \eqref{4.1}, we have $\mathscr{H}(x,\eta,p) \leq C(|\eta|)(1\,+\,p)$ and hence
\[
\mathscr{H}^{m-1}\Big(x,\eta,\frac{1}{2}\big|P-\mathscr{V}(x,\eta)\big|^2\Big)\, \leq\, C(|\eta|) \left( 1\, +\, \big|P-\mathscr{V}(x,\eta)\big|^{2m-2}\right).
\]
Further, by \eqref{4.1}, \eqref{4.2} and the above, we easily obtain the bounds
\begin{align} 
\label{bound1}
\big|H_P(x,\eta,P)\big|\, & \leq\, C(|\eta|) \big( 1\, +\, |P|^{2m-1}\big), 
\\ 
\label{bound2}
\big|H_\eta(x,\eta,P)\big|\, &\leq\, C(|\eta|) \big( 1\, +\, |P|^{2m}\big).
\end{align}
By standard results (see e.g.\ \cite{D}), \eqref{bound1}-\eqref{bound2} imply that the functional is Gateaux differentiable on $W^{1,2m}_b(\Om,\R^N)$ and the lemma follows.          \qed}

\ms

Now we show that the weak solutions $u^m$ of the respective Euler-Lagrange equations actually are smooth solutions. This will imply that the formal calculations of the previous section in the derivation of \eqref{6.7} make rigorous sense.

\bl[$C^2$ regularity] \label{lemma3} Let $(u^m)_1^\infty$ be the sequence of minimisers of Lemma \ref{lemma3}. Then, each $u^m$ is a classical solution in $C^2(\Om,\R^N)$ of the Euler-Lagrange equation \eqref{6.27}, and hence of the expanded form \eqref{6.7} of the same equation.
\el

\BPL \ref{lemma3}. Fix $m\geq 2$ and let us drop the superscripts and denote $u^m$ by just $u$. The first step is to prove higher local integrability and then bound the difference quotients of $\D u$ in $L^2_{\text{loc}}(\Om,\R^N)$. Let us fix $q\in \N$ and $\zeta \in C^\infty_c(\Om)$ with $0\leq \zeta \leq 1$. By recalling \eqref{W(u)}, we set:
\beq
\phi(x)\, :=\, \zeta(x) \int_{\inf \Om}^x \zeta(t) \big|\mathscr{W}u(t) \big|^q \mathscr{W}u(t)\,dt, \quad x\in \Om.
\eeq 
Then,  $\phi \in W^{1,1}_c(\Om,\R^N)$ and
\[
\begin{split}
\D \phi(x)\, =\, \zeta^2(x) & \big| \mathscr{W}u(x) \big|^q \mathscr{W}u(x)
\, +\,  \D \zeta(x) \int_{\inf \Om}^x \zeta(t) \big|\mathscr{W}u(t) \big|^q \mathscr{W}u(t)\,dt,
\end{split}
\]
for a.e.\ $x\in\Om$. Suppose now that $q\leq 2m-1$. Then, since $\D u\in L^{2m}(\Om,\R^N)$, we have that $\phi \in W^{1,2m}_c(\Om,\R^N)$. By inserting the test function $\phi$ in the weak formulation of the system \eqref{compressed} (i.e.\ \eqref{6.27}) and by suppressing again the arguments for the sake of brevity, we have
\[
\begin{split}
& \int_\Om \Bigg\{\mathscr{H}^{m-1}  \mathscr{H}_p \mathscr{W}u \cdot \Bigg[ \zeta^2 |\mathscr{W}u|^q \mathscr{W}u \,  +\,  \D \zeta \int_{\inf \Om} \zeta  |\mathscr{W}u |^q (\mathscr{W}u  )  \Bigg]\Bigg\}\\
&  +\, \int_\Om \Bigg\{ \mathscr{H}^{m-1}\mathscr{H}_p\,\Big( \mathscr{H}_\eta - \mathscr{H}_p (\mathscr{W}u )^\top  \mathscr{V}_\eta(\cdot,u)\Big) \cdot  \left[ \zeta \int_{\inf \Om} \zeta |\mathscr{W}u  |^q \mathscr{W}u    \right]\Bigg\}\, =\, 0.
\end{split}
\]
By \eqref{4.1}-\eqref{4.2}, we have $\mathscr{H}_p\geq c_0$ and $2h\geq c_0 |\mathscr{W}u|^2$. By using the bounds \eqref{bound1}, \eqref{bound2} (where $H$ is given by \eqref{6.22}), that $\zeta\leq1$ and the elementary inequalities
\[
\begin{split}
\int_{\inf \Om}^x|f| \, &\leq \, \int_\Om|f|, \ \quad x\in \Om,\ \ f\in L^1(\Om),\\
 t^{2m-1}\,& \leq \, t^{2m}+1, \quad t\geq0,
\end{split}
\]
we have 
\[
\begin{split}
\int_\Om \zeta^2|\mathscr{W}u|^{2m+q}\,  \leq & \, C\left(\int_\Om  \zeta|\mathscr{W}u|^{q+1} \right) \Bigg\{ \int_\Om |\D \zeta|\, \Big( \mathscr{H}^{m-1}\mathscr{H}_p\, |\mathscr{W}u|\Big)\ + \\
& +\, \int_\Om  \zeta \left( \mathscr{H}^{m-1} \mathscr{H}_p\,\Big| \mathscr{H}_\eta - \mathscr{H}_p (\mathscr{W}u )^\top  \mathscr{V}_\eta(\cdot,u) \Big| \right)\Bigg\}
\end{split}
\]
which gives
\[
\begin{split}
\int_\Om \zeta^2|\mathscr{W}u|^{2m+q}\, \leq & \, C   \big(\|u\|_{L^\infty(\Om)}\big)  \left(\int_\Om  \zeta|\mathscr{W}u|^{q+1} \right) \centerdot\\
&\centerdot  \int_\Om \Big\{|\D \zeta|\Big(1+ |\mathscr{W}u|^{2m-1}\Big)\,+\, \zeta \Big(1+ |\mathscr{W}u|^{2m}\Big) \Big\}.
\end{split}
\]
Hence, we have obtained
\beq \label{6.31}
\int_\Om \zeta^2|\mathscr{W}u|^{2m+q}\, \leq\, C\big(\|u\|_{L^\infty(\Om)}\big)  \left(\int_\Om  \zeta|\mathscr{W}u|^{q+1} \right)  \int_\Om 1+ |\mathscr{W}u|^{2m}.
\eeq
In view \eqref{6.31}, by taking $q+1=2m$ we have $\mathscr{W}u \in L^{4m-1}_{\text{loc}}(\Om,\R^N)$. Hence, we can iterate and choose $q+1=4m-1$ to find that $\phi \in W^{1,4m-1}_c(\Om,\R^N)$ which makes it admissible and we can repeat the process. Hence, by applying the estimate again we infer that $\mathscr{W}u \in L^{6m-2}_{\text{loc}}(\Om,\R^N)$.  By induction, the estimate holds for all integers of the form
\[
q\,=\, (2m-1)k,\quad k\in \N
\]
and we obtain that $\mathscr{W}u \in \bigcap_{r=1}^\infty L^r_{\text{loc}}(\Om,\R^N)$. In view of \eqref{W(u)} and since $u \in C^0(\overline{\Om},\R^N)$, we conclude that \[
\D u \,\in\, \bigcap_{r=1}^\infty L^r_{\text{loc}}(\Om,\R^N). 
\]
The next step is to prove that $\D^{1,t}\D u$ is bounded in $L^2_{\text{loc}}$. The idea is classical, but we provide the arguments for the sake of completeness. To this end, we test in the weak formulation of \eqref{compressed} against difference quotients of the form
\[
\phi\, :=\, -\D^{1,-t}\left(\zeta^2 \D^{1,t}u\right),\ \ \ \zeta \in C^\infty_c(\Om), \ \  \D^{1,t}u(x)\, =\, \frac{u(x+t)-u(x)}{t}, \ \ t\neq 0.
\]
Let $H$ be given by \eqref{6.22}. Then, for $\e>0$ and $t$ small, we have 
\beq \label{A}
\begin{split}
I\,:=\, &\left|\int_\Om \D^{1,t}\big(H_P(\cdot,u,\D u) \big)\cdot \Big(\zeta^2\D^{1,t}\D u\,+\,2\zeta \D \zeta\D^{1,t}u \Big) \right|\\
& \leq\, \int_\Om \big|H_\eta (\cdot,u,\D u) \big| \Big(\zeta |\D^{1,t}u|\, +\, \zeta^2|\D^{1,t}\D^{1,t}u| \Big)\\
&\leq K\left( \int_\Om \zeta \big|H_\eta (\cdot,u,\D u) \big|^2 \, +\, \int_\Om \zeta|\D u|^2   \right)\, +\, \e\int_\Om \zeta^2|\D^{1,t}\D u|^2. 
\end{split}
\eeq
for some constant $K>0$ independent of $t$. By using the inequality $H_{PP}\geq c_0I$ and the identity
\[
\begin{split}
\D^{1,t}&\big(H_P(\cdot,u,\D u) \big)(x)\\
 =&\, \int_0^1\Bigg\{H_{PP}\Big(\cdot,\la u(x+t)+(1-\la)u(x), \la \D u(x+t)+(1-\la)\D u(x)\Big)\D^{1,t}\D u(x) 
 \end{split}
\]
\[
\begin{split}
&\ \ \ \quad +\, H_{P\eta} \Big(\cdot,\la u(x+t)+(1-\la)u(x), \la \D u(x+t)+(1-\la)\D u(x)\Big)\D^{1,t}u(x)\\
&\ \ \ \quad +\, H_{Px_{1,t}} \Big(\cdot,\la u(x+t)+(1-\la)u(x), \la \D u(x+t)+(1-\la)\D u(x)\Big) \Bigg\} \, d\la
\end{split}
\]
(where $H_{Px_{1,t}} $ denotes difference quotient with respect to the $x$ variable), we have
\beq  \label{B}
\begin{split}
I \, \geq\, \frac{1}{K}\int_\Om \zeta^2|\D^{1,t}\D u|^2 \,-\, C\big(\|u\|_{L^\infty(\Om)}\big)\int_\Om\zeta \big|\P(|\D u|) \big|
\end{split}
\eeq
where $K>0$ is a constant independent of $t$, whilst $\P$ is a polynomial expression and it is a consequence of \eqref{4.1}-\eqref{4.2}. Since $u \in C^0(\overline{\Om},\R^N)$ and $\D u\in L^r_{\text{loc}}(\Om,\R^N)$ for all $r\geq 1$, by \eqref{A} and  \eqref{B} we obtain that $u\in  W^{2,2}_{\text{loc}}(\Om,\R^N)$. Thus, the calculations in the derivation of the expanded form of the system make sense a.e.\ on $\Om$. Since $H_{PP}$ is a strictly positive matrix, by a standard bootstrap argument in the system we obtain that $u \in C^2(\Om,\R^N)$ and the lemma follows.                \qed

\ms

Now we may prove the remaining assertions of our main result.

\ms

\noi \textbf{Proof of items \eqref{(2)}-\eqref{(4)} of Theorem \ref{theorem1}.} In view of Lemmas \ref{lemma1}, \ref{lemma2}, \ref{lemma3},  let $(u^m)_1^\infty$ denote the sequence of  minimisers in $C^0(\overline{\Om},\R^N) \cap C^2(\Om,\R^N)$ of the functionals \eqref{1.19} over the spaces $W^{1,2m}_b(\Om,\R^N)$. Then, along a subsequence
\beq \label{6.33}
\left\{
\begin{array}{l}
\ \ \  u^m \, -\!\!\!\!\larrow  u^\infty,\ \ \ \, \text{ in }C^0(\overline{\Om},\R^N),\ms\\
\D u^m \weak \D u^\infty,\ \ \text{ in } L^q(\Om,\R^N), \text{ for all }q\geq 1,
\end{array}
\right.
\eeq
as $m \ri \infty$, and the limit satisfies $u^\infty \in W^{1,\infty}_b(\Om,\R^N)$. Moreover, each $u^m$ is a classical solution of the system \eqref{6.7}, or equivalently of \eqref{3.23} with $f^\infty,F^\infty,A^\infty$ given by \eqref{3.21}-\eqref{3.22a}. The goal is to show that the limit map $u^\infty$ is a $\mD$-solution of the system \eqref{1.5a} with $\mF_\infty$ given by \eqref{1.14} (or equivalently \eqref{3.24}) and also that $u^\infty=b$ on $\p\Om$. We begin by observing that the boundary condition is satisfied as a result of the uniform convergence on $\overline{\Om}$. Moreover, by recalling \eqref{W(u)} and by multiplying \eqref{3.23} with $\D \big(\mathscr{W}{u^m} \big)$, we obtain 
\beq 
\label{4.21}
\begin{split}
  & \left\{\frac{ A^\infty(\cdot,{u^m},\D u^m)}{m-1}  +  \mathscr{H}_p^2\Big(\cdot,u,\frac{1}{2}\big|\mathscr{W}{u^m} \big|^2\Big)  \big|\mathscr{W}{u^m}\big|^2 I\right\} :   \D \big( \mathscr{W}{u^m} \big)
\\
 &  \ot  \D \big(\mathscr{W}{u^m}\big)\, = \, \left(\frac{f^\infty(\cdot,{u^m},\D u^m)}{m-1}+F^\infty(\cdot,{u^m},\D u^m) \right)\cdot  \D \big(\mathscr{W}{u^m}\big)
 \\
  & \hspace{55pt} \leq \, \left|\frac{f^\infty(\cdot,{u^m},\D u^m)}{m-1}+ F^\infty(\cdot,{u^m},\D u^m) \right|\big| \D \big(\mathscr{W}{u^m}\big)  \big|.
  \end{split}
\eeq
By \eqref{4.1} we have $\mathscr{H}_p \geq c_0$. In addition, by \eqref{3.22a} the matrix map $A^\infty$ is non-negative. Hence \eqref{4.21} gives the estimate
\beq \label{4.23}
\Big| \big| \mathscr{W}{u^m} \big|^2  \D \big(\mathscr{W}{u^m}\big) \Big| \, \leq\, \frac{1}{c_0^2}  \left|\frac{f^\infty(\cdot,{u^m},\D u^m)}{m-1}+F^\infty(\cdot,{u^m},\D u^m) \right|.
\eeq
By using the elementary inequality
\[
\left| \D \left( |f|^3\right) \right| \, \leq\, 3 \left| |f|^2 \D f\right|, \ \ \ f \in C^1(\Om,\R^N),
\]
\eqref{4.23} gives the estimate
\beq \label{4.24}
\Big| \D \left( | \mathscr{W}{u^m} |^3\right) \Big| \, \leq\, \frac{3}{c_0^2} \left|\frac{f^\infty(\cdot,{u^m},\D u^m)}{m-1}+F^\infty(\cdot,{u^m},\D u^m) \right|.
\eeq
By \eqref{4.24}, \eqref{6.33} and the form of the right hand side given by \eqref{3.21}, \eqref{3.22}, we have that the sequence
\beq \label{6.36a}
v^m\, :=\, | \mathscr{W}{u^m} |^3 \,=\, \big|  \D u^m -\mathscr{V}(\cdot,  u^m)  \big|^3
\eeq
is bounded in $W^{1,q}(\Om)$, for any $q\geq 1$. Hence, by the compactness of the imbedding $W^{1,q}(\Om) \Subset C^0(\overline{\Om})$, there is a continuous non-negative function $v^\infty$ such that 
\[
v^m \larrow v^\infty, \ \ \ \text{ in }C^0(\overline{\Om}),
\]
along perhaps a further subsequence as $m\ri \infty$. We claim that
\beq \label{6.37A}
| \mathscr{W}{u^\infty} |^3\,=\, \big|  \D u^\infty -\mathscr{V}(\cdot,  u^\infty)  \big|^3  \leq \, v^\infty, \ \ \text{a.e.\ on }\Om.
\eeq
Indeed, by \eqref{6.33} and the weak lower semi-continuity of the $L^3$ norm, for every $x\in \Om$ and $r>0$ fixed we have that
\beq \label{6.A}
\begin{split}
\frac{1}{2r}\int_{x-r}^{x+r}
\big| \mathscr{W}{u^\infty}  \big|^3 \, &\leq \ \underset{m\ri \infty}{\lim \inf} \,  \frac{1}{2r}\int_{x-r}^{x+r}  \big| \mathscr{W}{u^m} \big|^3\\
&= \ \lim_{m\ri \infty}\, \frac{1}{2r}\int_{x-r}^{x+r} v^m\\
&= \ \frac{1}{2r}\int_{x-r}^{x+r} v^\infty .
\end{split}
\eeq
By passing to the limit as $r\ri0$ in \eqref{6.A}, the Lebesgue differentiation theorem implies that the inequality  \eqref{6.37A} is valid a.e.\ on $\Om$. We now set
\[
\Om^\infty\, :=\, \big\{\, x\in \Om\ :\  v^\infty(x)>0\, \big\}.
\]
By the continuity of $v^\infty$, $\Om^\infty$ is open in $\Om$, the set $\Om\set\Om^\infty$ is closed in $\Om$ and 
\[
\Om\set\Om^\infty\, =\, \big\{\, x\in \Om\ :\  v^\infty(x)=0\, \big\}. 
\]
By \eqref{6.37A}, we have
\beq \label{6.37}
 \big| \mathscr{W}{u^\infty} \big| \, =\, 0, \ \ \text{a.e.\ on }\Om\set\Om^\infty.
\eeq
On the other hand, since $v^m \larrow v^\infty$ in $C^0(\overline{\Om})$, for any $U \Subset \Om^\infty$, there is a $\si_0>0$ and an $m(U)\in \N$ such that for all $m\geq m(U)$, {we have $v^m \geq \si_0$ on $U$} and hence by \eqref{6.36a}
\beq  \label{6.38}
| \mathscr{W}{u^m} | \, \geq\, (\si_0)^{\frac{2}{3}}, \ \  \text{ on }U.
\eeq
By \eqref{6.38}, \eqref{6.A} and \eqref{4.24}, we have
\beq  \label{6.39}
\Big|  \D \big( \mathscr{W}{u^m} \big) \Big| \, \leq\, \frac{3}{(c_0)^2  (\si_0)^{\frac{2}{3}}} \left|\frac{f^\infty(\cdot,{u^m},\D u^m)}{m-1}+F^\infty(\cdot,{u^m},\D u^m) \right|, \ \ \ \text{ on }\Om'.
\eeq
{By \eqref{6.39} and \eqref{6.33} we have} that $\D^2 u^m$ is bounded in $L^q_{\text{loc}}(\Om^\infty,\R^N)$. Hence, we have that
\[
\left\{
\begin{array}{l}
\ \ \ u^m \, \larrow  u^\infty,\ \ \ \ \text{ in }C^0(\Om^\infty,\R^N),\ms\\
{\, \D u^m \, \larrow \, \D u^\infty,\ \, \text{ in } L^q_{\text{loc}}(\Om^\infty,\R^N), \text{ for all }q\geq 1,} 
\ms\\
\D^2 u^m \! \weak \D^2 u^\infty, \text{ in } L^q_{\text{loc}}(\Om^\infty,\R^N), \text{ for all }q\geq 1. \ms
\end{array}
\right.
\]
Thus, by passing to the limit in the ODE system \eqref{6.7} as $m\ri \infty$ along a subsequence, we have that the restriction of $u^\infty$ over the open set $\Om^\infty$ is a strong a.e.\ solution of \eqref{3.24} on $\Om^\infty$. By bootstrapping in the equation, we have that actually $u^\infty\in C^2(\Om^\infty,\R^N)$. On the other hand, we have that 
\[
| \mathscr{W}  u^\infty |\, =\, 0, \quad \text{ a.e. on }\Om\setminus \Om^\infty. 
\]
Hence, if the set $\Om\set \Om^\infty$ has non-trivial topological interior, by differentiating the relation $ \D u^\infty = \mathscr{V}(\cdot,  u^\infty)$ we have that $\D^2 u^\infty$ exists a.e.\ on the interior of the open set $\Om\set \Om^\infty$ and by bootstrapping again we see that $u^\infty \in C^2 \big(\inter(\Om\set\Om^\infty),\R^N\big)$. Putting the above together, we have that $\D^2 u^\infty$ exists and is continuous on the open set $\Om_\infty$ defined in the statement of the theorem which is the union of $\Om^\infty$ and of the interior of $\Om\set\Om^\infty$:
\[
u^\infty\in C^2(\Om_\infty,\R^N),\quad \Om_\infty\, =\, \Om^\infty \cup \inter\, (\Om\set \Om^\infty).
\]
We now show that $u^\infty$ is a $\mD$-solution of \eqref{3.24} on $\Om$ (Definitions \ref{definition7}-\ref{definition11}). Let $\D^{1,h_i} \D u^\infty$ be the first difference quotients of $\D u^\infty$ along a sequence $h_i  \ri 0$ as $i\ri \infty$ and let $\mD^2{u^\infty}$ be a diffuse 2nd derivative of $u^\infty$ arising  from the subsequential weak* convergence of the  difference quotients, that is
\[
\de_{ \D^{1,h_{i_j} } \D u^\infty }  \weakstar \, \mD^2u^\infty, \ \ \ \text{ in } \mY \big(\Om, \smash{\overline{\R}}^N\big),
\]
as $j \ri \infty$, in the space of Young measures from $\Om\sub \R$ into the $1$-point compactification $\smash{\overline{\R}}^N=\R^N\cup\{\infty\}$. By the regularity of $u^\infty$ on $\Om^\infty$ and Lemma \ref{lemma10}, the restriction of any diffuse 2nd derivative on $\Om^\infty$ is the Dirac mass at the second derivatives: 
\beq \label{eq}
\mD^2u^\infty(x) \, =\, \de_{ \D^2 u^\infty(x)}, \ \ \ \text{ for a.e. } x\in \Om^\infty.
\eeq
Hence, $u^\infty$ is $\mD$-solution on $\Om^\infty$, since it is a strong solution on this subdomain. Consequently, for a.e.\ $x\in \Om^\infty \sub \Om$ and any $X \in \supp_*\big(\mD^2{u^\infty}(x)\big)$ we have
\beq
\begin{split} \label{4.31}
 &  \mathscr{H}_p^2   \Big(x, u^\infty(x),   \frac{1}{2}\big| \mathscr{W}u^\infty  (x)\big|^2\Big)  \big|\mathscr{W}u^\infty(x) \big|^2  \Big[\, X -\D \big(\mathscr{V}(\cdot,{u^\infty})\big)(x) \Big]  \\
& \, = \,  F^\infty\Big(x,{u^\infty}(x),\D u^\infty(x)\Big) .
  \end{split}
\eeq
Thus, $u^\infty$ is a $\mD$-solution of \eqref{1.5a} with $\mF_\infty$ given by \eqref{1.14} (i.e.\ \eqref{3.24}). On the other hand, since  $\big| \mathscr{W}u^\infty \big|=0$, a.e.\ on $\Om\set \Om^\infty$, for a.e.\ $x\in \Om\set \Om^\infty$ and any $X \in \supp_*\big(\mD^2{u^\infty}(x)\big)$ we have
 we have
\beq \label{4.32}
\begin{split}
&  \mathscr{H}_p^2   \Big(x, u^\infty(x),   \frac{1}{2}   \big| \mathscr{W}u^\infty   (x)\big|^2\Big)  \big| \mathscr{W}u^\infty(x) \big|^2   \Big[X-\D \big(\mathscr{W}(\cdot,{u^\infty}\big)(x) \Big]  \, =\, 0. 
  \end{split}
\eeq
Also, by \eqref{3.21} we see that the right hand side of \eqref{3.24} essentially vanishes on $\Om\set \Om^\infty$ as well: 
\beq \label{4.33}
F^\infty\big(\cdot,{u^\infty},\D u^\infty\big)\, =\, 0, \quad\text{ a.e. on }\Om\setminus \Om^\infty.
\eeq
By putting \eqref{4.31}, \eqref{4.32}, \eqref{4.33} together, we conclude that $u^\infty$ is indeed a $\mD$-solution of the Dirichlet problem for the fundamental equations in $L^\infty$, which is also a weak sequential limit of minimisers of the respective $L^m$ functionals as $m\ri \infty$ in the $W^{1,q}$ topology for any $q\geq 1$. In order to conclude it remains to establish the strong convergence of the derivatives of the sequence of minimisers $u^m$ to $u^\infty$. On the open set $\Om^\infty$ we have $Du^m \larrow Du^\infty$ in $L^q_{\text{loc}}(\Om^\infty,\R^N)$ and hence up to a further subsequence we have $Du^m(x) \larrow Du^\infty(x)$ for a.e.\ $x\in \Om^\infty$ as $m\ri\infty$. On the closed set $\Om \set \Om^\infty$, we have
\[
\begin{split}
\int_{\Om \set \Om^\infty} \big|\D u^m-\D u^\infty\big|^3\, \leq \, 9 \Bigg( & \int_{\Om \set \Om^\infty} \big|\D u^m- \mathscr{V}(\cdot,u^m)\big|^3 
\\
&+\, \int_{\Om \set \Om^\infty} \big|  \mathscr{V}(\cdot,u^m)- \mathscr{V}(\cdot,u^\infty)\big|^3 
\\
&+\, \int_{\Om \set \Om^\infty} \big|\mathscr{V}(\cdot,u^\infty) - \D u^\infty\big|^3 \Bigg),
\end{split}
\]
for any $q\in \N$. Since 
\[
\big|\D u^m- \mathscr{V}(\cdot,u^m)\big|^3\,=\, v^m \, \larrow \, v^\infty \, =\, \big|\D u^m- \mathscr{V}(\cdot,u^m)\big|^3\, =\, 0
\]
as $m\ri \infty$ in $C^0\big(\Om \set \Om^\infty,\R^N\big)$ and also $u^m \larrow u^\infty$ in $C^0\big(\overline{\Om},\R^N \big)$, we have that $\D u^m \larrow \D u^\infty$ in $L^3(\Om \set \Om^\infty,\R^N)$ along a subsequence as $m\ri \infty$. Conclusively, 
\[
Du^m(x) \larrow Du^\infty(x),\  \text{ for a.e.\ $x\in \Om$ as $m\ri\infty$ along a sequence} 
\]
and also 
\[
\| \D u^m \|_{L^q(\Om)}\, \leq\, C(q), \ \ \ \ q\in \N.
\]
Hence, if $E\sub \Om$ is measurable, we have the equi-integrability estimate
\[
\| \D u^m \|_{L^q(E)}\, \leq\,  \| \D u^m \|_{L^{q+1}(E)} |E|^{\frac{1}{q(q+1)}} \, \leq\, C(q+1) |E|^{\frac{1}{q(q+1)}}.
\]
The conclusion of strong convergence now follows from the above and the Vitali convergence theorem. Theorem \ref{theorem1} has been established.                 \qed

\ms

\ms

\noi \textbf{Acknowledgement.} The author is indebted to Jochen Br\"ocker for the scientific discussions on the subject of variational Data Assimilation.

\ms

\end{document}